\documentclass[10pt,english,french]{amsart}
\usepackage{helvet}
\usepackage{babel}
\usepackage[T1]{fontenc}
\usepackage[latin1]{inputenc}
\usepackage{a4wide}
\usepackage{amsmath,amsfonts,latexsym,amssymb}
\usepackage{xypic,pstricks,pst-node,pst-tree}
\xyoption{all}

 \theoremstyle{plain}    
 \newtheorem{thm}{Théorème}[section]
 \numberwithin{equation}{section} %% Comment out for sequentially-numbered
 \numberwithin{figure}{section} %% Comment out for sequentially-numbered
 \theoremstyle{plain}    
 \newtheorem*{thm*}{Théorème} 
 \theoremstyle{plain}    
 \newtheorem{prop}[thm]{Proposition} %%Delete [thm] to re-start numbering
 \theoremstyle{remark}    
 \newtheorem{notation}[thm]{Notation} 
 \theoremstyle{remark}
 \newtheorem{rem}[thm]{Remarque}
 \theoremstyle{plain}    
 \newtheorem{lem}[thm]{Lemme} %%Delete [thm] to re-start numbering

\theoremstyle{definition}
\newtheorem{enavant}[thm]{}

\begin{document}

\title{Autour d'une surface rationnelle dans $\mathbb{C}^{3}$ }

\maketitle
 \insert\footins{\footnotesize\noindent \textbf{Mathematics Subject Classification (2000)}: 14R10,14R25.


\noindent \textbf{Mots clefs}:  Surfaces de Danielewski, actions du  groupe additif, automorphismes de $\mathbb{C}^3$. \unskip\strut\par}

\begin{center}\author{{\bf A. DUBOULOZ} \\ \vspace{0.5cm}   \address Institut Fourier, Laboratoire de Mathématiques\\ UMR5582 (UJF-CNRS)\\ BP 74, 38402 ST MARTIN D'HERES CEDEX\\ FRANCE \\ \email adrien.dubouloz@ujf-grenoble.fr \\ \vspace{0.8cm}} \end{center}

\selectlanguage{english}
\begin{abstract}
Affine surfaces in $\mathbb{C}^{3}$ defined by an equation of the
form $x^{n}z-Q\left(x,y\right)=0$ have been increasingly studied
during the past 15 years. Of particular interest is the fact that
they come equipped with an action of the additive group $\mathbb{C}_{+}$
induced by such an action on the ambient space. The litterature of
the last decade may lead one to believe that there are essentially
no other of rational surfaces in $\mathbb{C}^{3}$ with this property.
In this note, we construct an explicit example of a surface nonisomorphic
to a one of the above type but equipped with a free $\mathbb{C}_{+}$-action
induced by an action on $\mathbb{C}^{3}$. We give an elementary and
self-contained proof of this fact. As an application, we construct
a wild but stably-tame automorphisme of $\mathbb{C}^{3}$ which seems
to be new. 
\end{abstract}
\selectlanguage{french}
\noindent

\begin{abstract}
Les surfaces affines de $\mathbb{C}^{3}$ définies par une équation
de la forme $x^{n}z-Q\left(x,y\right)=0$ ont fait l'objet d'une intense
activité de recherche durant la dernière décennie. Cette étude est
entre autres motivée par le fait que ces surfaces apparaissent de
manière parfois inattendue dans bon nombre de problèmes classiques
de la géométrie affine, tels que le Problème de Simplification, ou
les questions de linéarisation des actions de groupes algébriques
sur les espaces affines. Ainsi définies, ces surfaces se trouvent
être invariantes par certaines actions du groupe additif complexe
$\mathbb{C}_{+}$ sur $\mathbb{C}^{3}$. A la connaissance de l'auteur,
ces surfaces, ainsi que certaines de leurs généralisation naturelles,
sont essentiellement les seules possédant cette propriété à apparaître
dans la littérature. Afin de ne pas laisser s'installer le sentiment
que ces dernières sont effectivement les seules à jouir de cette propriété,
nous construisons un exemple explicite de surface de $\mathbb{C}^{3}$
non isomorphe à une du type précédent et invariante par une action
$\mathbb{C}_{+}$ sur $\mathbb{C}^{3}$. Nous en déduisons un automorphisme
polynomial sauvage de $\mathbb{C}^{3}$ apparemment inconnu jusqu'alors. 
\end{abstract}

\section*{Introduction}

La littérature de cette dernière décennie atteste que les surfaces
de $\mathbb{C}^{3}$ définies par une équation de la forme $x^{n}z-P\left(y\right)=0$,
ou plus généralement du type $x^{n}z-Q\left(x,y\right)=0$, ont bénéficié
d'une popularité croissante. L'intérêt porté à ces surfaces semble
remontrer à un désormais célèbre contre-exemple non publié au Problème
de Simplification de Zariski construit par W. Danielewski \cite{Dan89}.
Dans cet article, W. Danielewski établit en effet que les surfaces
$S_{1}$ et $S_{2}$ d'équations $xz-y\left(y-1\right)=0$ et $x^{2}z-y\left(y-1\right)=0$
ne sont pas isomorphes mais que $S_{1}\times\mathbb{C}$ est isomorphe
à $S_{2}\times\mathbb{C}$. Dans sa construction, W. Danielewski exploite
en particulier le fait que ces surfaces sont munies d'actions algébriques
du groupe additif complexe $\mathbb{C}_{+}$. Cela peut se voir directement
en utilisant la correspondance bien connue entre les actions de $\mathbb{C}_{+}$
sur une variété affine $S$ et les dérivations localement nilpotentes
de l'algèbre $\mathcal{O}\left(S\right)$ des fonctions régulières
sur $S$, une $\mathbb{C}$-dérivation $\partial$ de $\mathcal{O}\left(S\right)$
étant dite localement nilpotente si toute fonction régulière $f$
sur $S$ est annulée par une puissance convenable $n=n\left(f\right)$
de la dérivation. Dans le cas des surfaces $S_{1}$ et $S_{2}$ ci-dessus,
les actions considérées par W. Danielewski sont induites par les dérivations
localement nilpotentes $\partial_{1}=x\partial_{y}+2y\partial_{z}$
et $\partial_{2}=x^{2}\partial_{y}+2y\partial_{z}$ de $\mathbb{C}\left[x,y,z\right]$
qui annulent les équations de $S_{1}$et $S_{2}$ respectivement.
Du fait de la symétrie existante entre les variables $x$ et $z$
intervenant dans l'équation de $S_{1}$, cette dernière admet également
une seconde action de $\mathbb{C}_{+}$ dont les orbites générales
sont distinctes de celle de la première. Par contre, la surface $S_{2}$
ne jouit pas de cette propriété, L. Makar-Limanov \cite{ML01} ayant
en effet établi que les orbites générales de toute action de $\mathbb{C}_{+}$
sur une surface $S$ de la forme $x^{n}z-P\left(y\right)=0$, où $n\geq2$,
coïncident avec les fibres générales de la projection ${\rm pr}_{x}:S\rightarrow\mathbb{C}$.
La même année, L. Makar-Limanov et T. Bandman \cite{BML01} ont prouvé
qu'une surface admettant deux actions de $\mathbb{C}_{+}$-actions
à orbites générales distinctes et plongeable dans $\mathbb{C}^{3}$
est isomorphe à une surface définie par une équation du type $xz-P\left(y\right)=0$,
où $P$ est un polynôme non constant.

Depuis lors, les surfaces isomorphes à des surfaces du type $S_{Q,n}$
définies par une équation de la forme $x^{n}z-Q\left(x,y\right)=0$
dans $\mathbb{C}^{3}$ ont fait l'objet de nombreux travaux généralisant
les résultats obtenus par L. Makar-Limanov pour le cas des surfaces
du type $x^{n}z-P\left(y\right)=0$. L'étude à équivalence algébrique
près de leurs plongements dans l'espace affine $\mathbb{C}^{3}$ menée
successivement par G. Freudenburg et L. Mauser-Jauslin \cite{FrMo02},
L. Mauser-Jauslin et P.M. Poloni \cite{MoP05} a conduit en particulier
à de nouveaux résultats concernant la nature du groupe des automorphismes
de $\mathbb{C}^{3}$. Une des propriétés remarquables partagée par
toutes les surfaces du type $S_{Q,n}$ est d'admettre une action de
$\mathbb{C}_{+}$ induite par une action sur $\mathbb{C}^{3}$ pour
laquelle le morphisme quotient algébrique $q:S_{Q,n}\rightarrow S_{Q,n}/\!/\mathbb{C}_{+}$
s'identifie à la restriction sur $S_{Q,n}$ de la projection ${\rm pr}_{x}:\mathbb{C}^{3}\rightarrow\mathbb{C}$.
Lorsque le polynôme $Q\left(0,y\right)$ n'est pas constant, la fibration
quotient $q$ est surjective et se restreint en un $\mathbb{A}^{1}$-fibré
trivial au-dessus de $\mathbb{C}^{*}$. Il est alors naturel de chercher
à savoir si les surfaces $S_{Q,n}$ sont les seules surfaces rationnelles
dans $\mathbb{C}^{3}$ possédant cette propriété. A la connaissance
de l'auteur, il n'est pas fait mention dans la littérature récente
d'un exemple venant infirmer cette conjecture raisonnable. Le propos
de cet article est par conséquent de présenter un contre-exemple explicite,
à savoir le suivant:

\begin{thm*}
La surface $S$ de $\mathbb{C}^{3}=\textrm{Spec}\left(\mathbb{C}\left[x,y,z\right]\right)$
définie par l'équation \[
x\left(xz^{2}-2y^{2}z+5z\right)+y^{4}-5y^{2}+4=0\]
 n'est pas isomorphe à une surface du type $S_{Q,n}$. Elle est néanmoins
invariante par l'action du groupe additif $\mathbb{C}_{+}$ sur $\mathbb{C}^{3}$
définie par la dérivation localement nilpotente \[
\partial=x\left(2xz-2y^{2}+5\right)\frac{\partial}{\partial y}+\left(x+2y\left(2xz-2y^{2}+5\right)\right)\frac{\partial}{\partial z}.\]
 De plus, la fibration quotient algébrique $q:S\rightarrow S/\!/\mathbb{C}_{+}$
associée à la $\mathbb{C}_{+}$-action induite sur $S$ coïncide avec
la projection ${\rm pr}_{x}\mid_{S}:S\rightarrow\mathbb{C}$ et se
restreint en un $\mathbb{A}^{1}$-fibré trivial au-dessus de $\mathbb{C}^{*}$. 
\end{thm*}
\indent\newline\indent Les spécialistes se convaincront aisément
que ce résultat pourrait se déduire presque immédiatement de certains
faits généraux sur les surfaces avec $\mathbb{C}_{+}$-actions établis
par l'auteur \cite{Dub02,DubG03} et P.M. Poloni \cite{DubP04}, faits
qui ne sont pour l'essentiel que des généralisations ou des reformulations
de résultats antérieurs dus à J. Bertin \cite{Ber83}, L. Makar-Limanov
\cite{ML01} et K.H. Fieseler \cite{Fie94}. Notre démarche ici est
toute autre puisqu'elle consiste à donner une preuve complète directe
et essentiellement élémentaire du résultat ci-dessus, preuve dont
le but principal est d'illustrer sur un exemple particulier l'usage
de ces méthodes générales. 

En guise d'application, nous construisons ensuite un automorphisme
de $\mathbb{C}^{3}$, qui d'après un critère récemment établi par
U. Umirbaev et I. Shestakov \cite{UmS04-2} se trouve être un automorphisme
\emph{sauvage}, i.e. un automorphisme n'admettant pas de décomposition
en produit d'automorphismes affines et triangulaires de $\mathbb{C}^{3}$.

\section{Une surface rationnelle $\mathbb{C}_{+}$-invariante dans $\mathbb{C}^{3}$}

\noindent L'objectif de cette section est de donner une preuve directe
complète du résultat suivant:

\begin{thm}
\label{thm:Th_Principal} La surface $S$ de $\mathbb{C}^{3}=\textrm{Spec}\left(\mathbb{C}\left[x,y,z\right]\right)$
définie par l'équation \[
x\left(xz^{2}-2y^{2}z+5z\right)+y^{4}-5y^{2}+4=0\]
 n'est pas isomorphe à une surface du type $S_{Q,n}$. Elle est néanmoins
invariante par l'action du groupe additif $\mathbb{C}_{+}$ sur $\mathbb{C}^{3}$
définie par la dérivation localement nilpotente \[
\partial=x\left(2xz-2y^{2}+5\right)\frac{\partial}{\partial y}+\left(x+2y\left(2xz-2y^{2}+5\right)\right)\frac{\partial}{\partial z}.\]
 
\end{thm}
\begin{enavant} Le lecteur se convaincra sans peine que la dérivation
$\partial$ de $\mathbb{C}\left[x,y,z\right]$ proposée ci-dessus
est bien localement nilpotente et qu'elle annule l'équation de la
surface $S$. Par conséquent, $S$ est bien invariante par l'action
de $\mathbb{C}_{+}$ correspondante, comme annoncé. D'autre part,
on observe que l'équation proposée ne peut être mise sous une forme
du type $x^{n}z-Q\left(x,y\right)=0$ dans le système de coordonnées
$\left(x,y,z\right)$ utilisé. Rien ne garantit cependant que $S$
ne soit pas simplement un plongement ``exotique'' dans $\mathbb{C}^{3}$
d'une surface $S_{Q,n}$ de ce type. 

\end{enavant}

\begin{enavant} Il s'agit par conséquent de montrer que $S$, considérée
comme une surface abstraite, n'est pas isomorphe à une du type $S_{Q,n}$.
Bien qu'il soit certainement possible d'établir ce résultat par d'autres
arguments, notre approche repose sur l'étude de la structure des $\mathbb{A}^{1}$-\emph{fibrations}
sur ces différentes surfaces, c'est-à-dire des morphismes surjectifs
$q:S_{Q,n}\rightarrow\mathbb{C}$ dont la fibre générique est isomorphe
à une droite affine. Cela se justifie entre autres par le fait que
cette méthode, qui peut se généraliser agréablement à des situations
de départ plus complexes, a le mérite d'illustrer un certain nombre
de techniques à l'origine de progrès récents dans l'étude des surfaces
affines rationnelles. Le principe de la démonstration est le suivant
: on montre tout d'abord que toute $\mathbb{A}^{1}$-fibration sur
$S$ s'identifie à automorphismes de la base près au morphisme quotient
algébrique $\pi:S\rightarrow S/\!/\mathbb{C}_{+}$ associé à l'action
de $\mathbb{C}_{+}$ déduite de la dérivation localement nilpotente
$\partial$ ci-dessus. Sachant cela, on établit qu'à un changement
de base près, un isomorphisme $\Phi:S\stackrel{\sim}{\rightarrow}S_{Q,n}$
entre $S$ et une surface de type $S_{Q,n}$ munie de sa $\mathbb{A}^{1}$-fibration
${\rm pr}_{x}:S_{Q,n}\rightarrow\mathbb{A}^{1}$ doit être un isomorphisme
de surfaces fibrées. Il ne reste alors plus qu'à vérifier qu'aucune
$\mathbb{A}^{1}$-fibration ${\rm pr}_{x}:S_{Q,n}\rightarrow\mathbb{A}^{1}$
sur une surface du type $S_{Q,n}$ n'a la même structure que la fibration
quotient $\pi:S\rightarrow S/\!/\mathbb{C}_{+}$. 

\end{enavant}

\begin{enavant} Ce plan étant établi, il suffirait pour démonter
le Théorème \ref{thm:Th_Principal} de remarquer que $S$ est une
surface de Danielewski au sens de \cite{DubG03} par rapport au morphisme
quotient $\pi:S\rightarrow S/\!/\mathbb{C}_{+}$. En combinant alors
la classification de ces surfaces donnée dans \emph{loc. cit.} avec
le fait que toute surface admettant deux $\mathbb{A}^{1}$-fibrations
à fibres générales distinctes et plongeable dans $\mathbb{C}^{3}$
est isomorphe à une surface d'équation $xz-P\left(y\right)=0$ pour
un certain polynôme non constant $P$ (cf. \cite{BML01}), on pourrait
conclure de manière directe que $\pi$ est l'unique $\mathbb{A}^{1}$-fibration
sur $S$ à automorphismes de la base près. Le fait que $S$ ne soit
pas isomorphe en tant que surface fibrée à une surface du type $S_{Q,n}$
se déduirait ensuite de la description des $\mathbb{A}^{1}$-fibrations
sur certaines de ces surfaces obtenue par P.M. Poloni et l'auteur
dans \cite{DubP04}. Néanmoins, afin d'illustrer dans un cas particulier
les différentes méthodes mises en oeuvres dans l'étude générale des
surfaces affines avec actions de groupes additifs, nous allons donner
une preuve complète du Théorème \ref{thm:Th_Principal}, s'inspirant
certes très fortement de ces résultats antérieurs, mais essentiellement
directe et indépendante. 

\end{enavant}

\subsection{Structures $\mathbb{A}^{1}$-fibrées sur $S$ }

\indent\newline\noindent Dans ce paragraphe, nous décrivons la structure
du morphisme quotient $\pi:S\rightarrow S/\!/\mathbb{C}_{+}$ sur
$S$ pour la $\mathbb{C}_{+}$-action induite par la dérivation $\partial$.
Nous montrons en particulier que cette action est libre et que le
``bon'' quotient à considérer ici s'identifie à une droite affine
à quatre origines. Commençons par introduire une $\mathbb{A}^{1}$-fibration
``naturelle'' sur $S$ qui, comme nous le verrons en \ref{pro:Quotient_algebrique}
ci-dessous, s'identifie au morphisme quotient algébrique $\pi:S\rightarrow S/\!/\mathbb{C}_{+}$. 

\begin{prop}
\label{pro:Fibration_sur_S} La surface $S\subset\mathbb{C}^{3}$
définie par l'équation \[
x\left(xz^{2}-2y^{2}z+5z\right)+y^{4}-5y^{2}+4=0\]
 est irréductible et non singulière. La projection ${\rm pr}_{x}:\mathbb{C}^{3}\rightarrow\mathbb{C}$
induit sur $S$ une $\mathbb{A}^{1}$-fibration $\pi:S\rightarrow\mathbb{C}$
satisfaisant les propriétés suivantes. 

1) $\pi^{-1}\left(\mathbb{C}^{*}\right)$ est isomorphe au $\mathbb{A}^{1}$-fibré
trivial $\mathbb{C}^{*}\times\mathbb{C}$ de base $\mathbb{C}^{*}$.
De plus, la restriction du polynôme $f=xz-y^{2}$ sur $S$ induit
une coordonnée sur le second facteur du produit. 

2) La fibre $\pi^{-1}\left(0\right)$ est réduite, constituée de la
réunion disjointe de quatre droites affines $C_{\pm1}$ et $C_{\pm2}$
d'équations respectives $\left\{ x=0\,;y=\pm1\right\} $ et $\left\{ x=0\,;y=\pm2\right\} $. 
\end{prop}
\begin{proof}
Il suffit d'établir que $\pi:S\rightarrow\mathbb{C}$ est une $\mathbb{A}^{1}$-fibration
satisfaisant les propriétés annoncées. En effet, on déduit de la forme
de l'équation que $x$ n'est pas un diviseur de zéro dans l'anneau
$B=\mathcal{O}\left(S\right)$. Par conséquent, l'homomorphisme canonique
de localisation $B\rightarrow B_{x}=B\otimes_{\mathbb{C}\left[x\right]}\mathbb{C}\left[x,x^{-1}\right]$
est injectif. Or 1) implique en particulier que $B_{x}$ est intègre,
si bien que $B$ le sera aussi. La lissité de $S$ provient quant
à elle du fait que la $\mathbb{A}^{1}$-fibration $\pi$ est un morphisme
lisse. La seconde assertion se déduit immédiatement de la ré-écriture
de l'équation définissant $S$ sous la forme $xy-\left(f+1\right)\left(f+4\right)=0$
. On observe de plus que \[
B_{x}\simeq\mathbb{C}\left[x,x^{-1}\right]\left[y,z\right]/\left(xy-\left(f+1\right)\left(f+4\right)\right)\simeq\mathbb{C}\left[x,x^{-1}\right]\left[f\right].\]
 Cette dernière algèbre est elle-même isomorphe à un anneau de polynômes
en une variable au-dessus de $\mathbb{C}\left[x,x^{-1}\right]$ puisque
$f$ et $x$ sont algébriquement indépendants. On conclut donc que
$\pi$ est bien une $\mathbb{A}^{1}$-fibration satisfaisant qui plus
est la propriété 1).
\end{proof}
\begin{notation}
Nous écrirons dorénavant l'équation définissant la surface $S$ sous
la forme $P=xy-\left(f+1\right)\left(f+4\right)=0$, $f$ désignant
le polynôme $xz-y^{2}$. 
\end{notation}

\subsubsection{Fibration ``quotient algébrique''}

\indent\newline\noindent Dans ce paragraphe, nous étudions la nature
de l'action de $\mathbb{C}_{+}$ sur $S$ induite par la dérivation
localement nilpotente $\partial$ de $\mathbb{C}\left[x,y,z\right]$
définie dans le Théorème \ref{thm:Th_Principal}.

\begin{enavant} Plutôt que de décrire immédiatement cette action
sur $S$, il est plus instructif d'étudier au préalable l'action correspondante
sur $\mathbb{C}^{3}$. Puisque $\partial$ est par définition une
$\mathbb{C}\left[x\right]$-dérivation, elle s'étend canoniquement
en une $\mathbb{C}\left[x,x^{-1}\right]$-dérivation $\tilde{\partial}$
de l'anneau $\mathbb{C}\left[x,x^{-1},y,z\right]$. En ré-exprimant
$\partial$ en terme du polynôme $f=xz-y^{2}$, on obtient pour $\tilde{\partial}$
la forme compacte suivante :\begin{equation}
\tilde{\partial}=x\left(2f+5\right)\frac{\partial}{\partial y}+\left(x+2y\left(2f+5\right)\right)\frac{\partial}{\partial z}=-\frac{\partial P}{\partial z}\frac{\partial}{\partial y}+\frac{\partial P}{\partial y}\frac{\partial}{\partial z}=x^{2}\frac{\partial}{\partial f}.\label{eq:partial_tilde_eqn}\end{equation}
 En effet, $f$ et $P$ sont des ``variables'' de $\mathbb{C}\left[x,x^{-1},y,z\right]$,
au sens que $\mathbb{C}\left[x,x^{-1}\right]\left[y,z\right]\simeq\mathbb{C}\left[x,x^{-1}\right]\left[f,P\right]$,
pour lequelles on a les égalités $\tilde{\partial}f=x^{2}$et $\tilde{\partial}P=0$
par construction. Ainsi, $\textrm{Ker}\left(\tilde{\partial}\right)=\mathbb{C}\left[x,x^{-1},P\right]$.
On en conclut aisément que $\textrm{Ker}\left(\partial\right)=\mathbb{C}\left[x,P\right]$
car $x$ ne divise pas $P$. 

\end{enavant}

\begin{enavant}\label{txt:Lieu_fixe_de_partial} L'algèbre $\mathbb{C}\left[x,y,z\right]^{\mathbb{C}_{+}}$
des fonctions sur $\mathbb{C}^{3}$ invariantes sous la $\mathbb{C}_{+}$-action
définie par $\partial$ s'identifiant au noyau de $\partial$, on
conclut que $\mathbb{C}\left[x,y,z\right]^{\mathbb{C}_{+}}=\mathbb{C}\left[x,P\right]$,
ce dernier étant par ailleurs isomorphe à un anneau de polynômes en
deux variables puisque $x$ et $P$ sont clairement algébriquement
indépendants. Le morphisme quotient $q:\mathbb{C}^{3}\rightarrow\mathbb{C}^{3}/\!/\mathbb{C}_{+}\simeq\mathbb{C}^{2}$
est donc défini par $\left(x,y,z\right)\mapsto\left(x,P\left(x,y,z\right)\right)$.
Il est en particulier surjectif puisque $x$ ne divise pas $P$. On
retrouve explicitement dans cet exemple deux résultats généraux portant
sur les actions de $\mathbb{C}_{+}$ sur $\mathbb{C}^{3}$, à savoir,
le fait que le quotient algébrique $\mathbb{C}^{3}/\!/\mathbb{C}_{+}$
soit isomorphe à $\mathbb{C}^{2}$, et la propriété de surjectivité
du morphisme quotient, établis respectivement par M. Miyanishi \cite{Miy85}
et P. Bonnet \cite{Bon02}. D'autre part, la relation (\ref{eq:partial_tilde_eqn})
ci-dessus implique que $\partial\left(\mathbb{C}\left[x,y,z\right]\right)\subset\left(x,2y^{2}-5\right)\mathbb{C}\left[x,y,z\right]$.
On en conclut que la $\mathbb{C}_{+}$-action correspondante n'est
pas libre et que ses points fixes sont exactement les points appartenant
à la réunion des droites d'équations $L_{\pm}=\left\{ x=0;\, y=\pm\sqrt{5/2}\right\} $. 

\end{enavant}

\begin{prop}
\label{pro:Quotient_algebrique} L'action du groupe $\mathbb{C}_{+}$
sur $S$ induite par la dérivation localement nilpotente $\partial$
est libre. Son quotient algébrique $\pi:S\rightarrow S/\!/\mathbb{C}_{+}$
s'identifie à la projection ${\rm pr}_{x}\mid_{S}:S\rightarrow\mathbb{C}$.
\end{prop}
\begin{proof}
Il est clair que les droites de points fixes $L_{+}$ et $L_{-}$
ne rencontrent pas $S$. L'action $\mathbb{C}_{+}$ induite est par
conséquent libre. D'autre par, puisque $\textrm{Ker}\left(\partial\right)=\mathbb{C}\left[x,y,z\right]^{\mathbb{C}_{+}}=\mathbb{C}\left[x,P\right]$,
on conclut que $\mathcal{O}\left(S\right)^{\mathbb{C}_{+}}=\left(\mathbb{C}\left[x,y,z\right]/\left(P\right)\right)^{\mathbb{C}_{+}}=\mathbb{C}\left[x\right]$.
Ainsi, le morphisme quotient $\pi$ est induit par l'inclusion $\mathbb{C}\left[x\right]\hookrightarrow\mathcal{O}\left(S\right)$,
ce qui termine la preuve. 
\end{proof}

\subsubsection{Fibration ``quotient géométrique''}

\indent\newline\noindent Puisque $\pi:S\rightarrow S/\!/\mathbb{C}_{+}\simeq\mathbb{C}$
s'identifie à la projection ${\rm pr}_{x}:S\rightarrow\mathbb{C}$
qui est elle-même une $\mathbb{A}^{1}$-fibration triviale au dessus
de $\mathbb{C}^{*}$ d'après la Proposition \ref{pro:Fibration_sur_S},
on conclut que la restriction $\pi:\pi^{-1}\left(\mathbb{C}^{*}\right)\rightarrow\mathbb{C}^{*}$
de $\pi$ au-dessus de $\mathbb{C}^{*}$ est un fibré principal homogène
sous $\mathbb{C}_{+}$. Par contre, cette structure de fibré principal
homogène sous $\mathbb{C}_{+}$ ne s'étend pas à $\pi:S\rightarrow\mathbb{C}$
puisque la fibre $\pi^{-1}\left(0\right)$ est constituée d'une réunion
disjointe de quatre droites affines complexes. Nous allons cependant
construire un schéma $X$ et un morphisme $\rho:S\rightarrow X$ ayant
la propriété d'avoir la structure d'un fibré principal homogène de
groupe $\mathbb{C}_{+}$ et de base $X$ pour l'action considérée.
Cette construction est basée sur une observation faite par W. Danielewski
\cite{Dan89} et développée ensuite par K. H. Fieseler \cite{Fie94}
puis par l'auteur \cite{DubG03}. 

\begin{enavant} Le schéma $X$ que nous recherchons est obtenu en
remplaçant l'origine de $\mathbb{C}$ par quatre points $x_{\pm1}$
et $x_{\pm2}$ correspondant respectivement aux composantes $C_{\pm1}$
et $C_{\pm4}$ de la fibre $\pi^{-1}\left(0\right)$. Autrement dit,
$X$ est construit en recollant quatre copies $\delta_{\alpha}:X_{\alpha}\stackrel{\sim}{\rightarrow}\mathbb{C}$
de $\mathbb{C}$, où $\alpha=\pm1,\,\pm2$, par l'identité en dehors
de leurs origines respectives $\delta_{\alpha}^{-1}\left(0\right)=x_{\alpha}$.
Ce schéma non séparé $X$, isomorphe à la droite affine à quatre origines,
est naturellement muni d'un morphisme $\delta:X\rightarrow\mathbb{C}$
induisant les isomorphismes $\delta_{\alpha}$ sur les ouverts du
recouvrement canonique de $X$ par les $X_{\alpha}$. Il existe alors
un unique morphisme $\rho:S\rightarrow X$ factorisant $\pi:S\rightarrow\mathbb{C}$
au sens que $\pi=\delta\circ\rho$ et tel que $\rho^{-1}\left(x_{\alpha}\right)=C_{\alpha}$
pour tout $\alpha=\pm1,\,\pm2$. Par construction, les fibres du morphisme
$\rho:S\rightarrow X$ coïncident avec les orbites de l'action de
$\mathbb{C}_{+}$ sur $S$ induit par la dérivation $\partial$. L'existence
sur $\rho:S\rightarrow X$ d'une structure de fibré principal homogène
sous $\mathbb{C}_{+}$ de base $X$ nous est alors normalement garantie
par un argument standard de descente (cf \cite{SGA1}, exposé XI).
Néanmoins, nous allons donner dans le paragraphe suivant une description
explicite complète de cette structure qui nous sera utile par la suite. 

\end{enavant}

\begin{rem}
La construction élémentaire précédente se justifie de la façon plus
conceptuelle suivante. L'action considérée étant libre, on sait a
priori qu'un vrai quotient $\rho:S\rightarrow S/\mathbb{C}_{+}=\mathcal{C}$,
au sens que $S$ devient via $\rho$ un fibré principal homogène sous
$\mathbb{C}_{+}$ de base $\mathcal{C}$, existe dans la catégorie
des espaces algébriques. Dans ce cas particulier, l'espace algébrique
$\mathcal{C}$ est non singulier de dimension $1$, donc est un schéma,
éventuellement non séparé. La propriété universelle du quotient $\rho$
implique que le morphisme quotient algébrique $\pi:S\rightarrow\mathbb{C}$
se factorise par un morphisme $\delta:\mathcal{C}\rightarrow\mathbb{C}$
induisant nécessairement un isomorphisme $\delta:\delta^{*}\left(\mathbb{C}^{*}\right)\stackrel{\sim}{\rightarrow}\mathbb{C}^{*}$.
En effet, on sait déjà que $\pi$ induit une structure de fibré principal
homogène sur $\pi^{-1}\left(\mathbb{C}^{*}\right)$. Puisque les fibres
du morphisme $\rho$ doivent être connexes, on conclut que le candidat
évident pour $\mathcal{C}$ est la droite affine à quatre origines
$X$ construite ci-dessus. 
\end{rem}

\subsubsection{Structure de fibré principal homogène sur $S$ }

\indent\newline\noindent Dans ce paragraphe, nous allons expliciter,
à des fins d'utilisation future, la structure $\rho:S\rightarrow X$
de fibré principal homogène sous $\mathbb{C}_{+}$. 

\begin{enavant} \label{txt:Trivilisations_locales} Par construction,
pour tout $\alpha=\pm1,\,\pm2$, la pré-image $\rho^{-1}\left(X_{\alpha}\right)$
de l'ouvert canonique $X_{\alpha}\simeq\mathbb{C}$ de $X$ est isomorphe
au complémentaire $S_{\alpha}=S\setminus{\displaystyle \bigcup_{\alpha'\neq\alpha}}C_{\alpha'}$
dans $S$ des composantes $C_{\alpha'}$ de la fibre $\pi^{-1}\left(0\right)$
distinctes de $C_{\alpha}$. D'après la Proposition \ref{pro:Fibration_sur_S},
le polynôme $f=xz-y^{2}$ induit un isomorphisme $\mathcal{O}\left(S\right)_{x}=\mathcal{O}\left(S\right)\otimes_{\mathbb{C}\left[x\right]}\mathbb{C}\left[x,x^{-1}\right]\simeq\mathbb{C}\left[x,x^{-1}\right]\left[f\right]$,
tandis que $f\mid_{C_{\alpha}}=-\alpha^{2}$ pour tout $\alpha=\pm1,\,\pm2$.
Par conséquent, la restriction sur $S_{\alpha}$ de la fonction rationnelle
$x^{-1}\left(f+\alpha^{2}\right)\in\mathbb{C}\left[x,x^{-1},y,z\right]$
est une fonction régulière $g_{\alpha}\in\mathcal{O}\left(S_{\alpha}\right)$.
De plus, puisque $S$ est définie par l'équation $xy-\left(f+1\right)\left(f+4\right)=0$,
on conclut que $g_{\alpha}$ s'identifie à la restriction sur $S_{\alpha}$
de l'une des fonctions rationnelles $\left(f+4\right)^{-1}y$ si $\alpha=\pm1$
et $\left(f+1\right)^{-1}y$ si $\alpha=\pm2$. On en déduit que $g_{\alpha}\mid_{C_{\alpha}}=\left(-1\right)^{\left|\alpha\right|-1}\alpha/3$.
Ainsi, pour tout $\alpha=\pm1,\,\pm2$, la restriction sur $S_{\alpha}$
de la fonction rationnelle \[
x^{-2}\left(f+\left(-1\right)^{\left|\alpha\right|}\frac{\alpha}{3}x+\alpha^{2}\right)\in\mathbb{C}\left[x,x^{-1},y,z\right]\]
 est une fonction régulière $h_{\alpha}\in\mathcal{O}\left(S_{\alpha}\right)$. 

\end{enavant}

\begin{lem}
Pour tout $\alpha=\pm1,\:\pm2$, l'homomorphisme de $\mathbb{C}\left[x\right]$-algèbres
$\tau_{\alpha}^{*}:\mathbb{C}\left[x\right]\left[u_{\alpha}\right]\rightarrow\mathcal{O}\left(S_{\alpha}\right)$,
$u_{\alpha}\rightarrow h_{\alpha}$ est un isomorphisme. 
\end{lem}
\begin{proof}
Par construction, $\tau_{\alpha}^{*}$ induit un isomorphisme $\mathbb{C}\left[x,x^{-1}\right]\left[u_{\alpha}\right]\stackrel{\sim}{\rightarrow}\mathcal{O}\left(S_{\alpha}\right)_{x}$
puisque $\mathcal{O}\left(S_{\alpha}\right)_{x}\simeq\mathcal{O}\left(S\right)_{x}\simeq\mathbb{C}\left[x,x^{-1}\right]\left[f\right]\simeq\mathbb{C}\left[x,x^{-1}\right]\left[h_{\alpha}\right]$
d'après la Proposition \ref{pro:Fibration_sur_S}. Il suffit donc
pour conclure de vérifier que $h_{\alpha}$ induit une fonction coordonnée
sur $C_{\alpha}\simeq\mathbb{C}$, ce qui sera le cas si $h_{\alpha}$
n'est pas constante sur $C_{\alpha}$. Pour cela, on peut se contenter,
étant donnée la définition de $h_{\alpha}$, de vérifier que la restriction
de la $2$-forme $dx\wedge df$ sur $S_{\alpha}$ s'annule exactement
à l'ordre $2$ le long de $C_{\alpha}$. Or on établit aisément que
les $2$-formes $x^{2}dx\wedge dz$ et $\left(x+2y\left(2f+5\right)\right)dx\wedge df$
induise la même forme $\omega\in\Omega^{2}\left(S\right)$. Or la
restriction sur $S$ du polynôme $x+2y\left(2f+5\right)$ ne s'annule
pas le long de la fibre $\pi^{-1}\left(0\right)$, tandis que $z$
induit une fonction coordonnée sur chaque composante de $\pi^{-1}\left(0\right)$
d'après la preuve de la Proposition \ref{pro:Fibration_sur_S}. On
en conclut que $\textrm{ord}_{C_{\alpha}}\left(dx\wedge df\right)=\textrm{ord}_{C_{\alpha}}\left(x^{2}dx\wedge dz\right)=2$,
ce qui termine la preuve. 
\end{proof}
\begin{enavant} \label{txt:Fonctions_de_transition} Les trivialisations
$\tau_{\alpha}^{*}:\mathbb{C}\left[x\right]\left[u_{\alpha}\right]\stackrel{\sim}{\rightarrow}\mathcal{O}\left(S_{\alpha}\right)$,
$\alpha=\pm1,\,\pm2$ construites ci-dessus permettent d'identifier
la surface $S$ au schéma ${\displaystyle \bigsqcup_{\alpha=\pm1,\pm2}X_{\alpha}\times\mathbb{C}/\sim}$,
obtenu en recollant quatre copies $X_{\alpha}\times\mathbb{C}$ de
$\mathbb{C}^{2}$ le long de leurs ouverts respectifs $\left(X_{\alpha}\setminus\left\{ 0\right\} \right)\times\mathbb{C}$
via les isomorphismes \[
\left(X_{\alpha}\setminus\left\{ 0\right\} \right)\times\mathbb{C}\ni\left(x,u_{\alpha}\right)\rightarrow\left(x,u_{\alpha}+g_{\alpha\alpha'}\right)\in\left(X_{\alpha'}\setminus\left\{ 0\right\} \right)\times\mathbb{C},\]
 où, pour tout couple $\left(\alpha,\alpha'\right)\in\left\{ \pm1,\,\pm2\right\} $
on a posé \[
g_{\alpha\alpha'}=x^{-2}\left(\left(\left(\alpha'\right)^{2}-\alpha^{2}\right)+\frac{x}{3}\left(\left(-1\right)^{\left|\alpha'\right|}\alpha'-\left(-1\right)^{\left|\alpha\right|}\alpha\right)\right)\in\mathbb{C}\left[x,x^{-1}\right].\]
 En identifiant à leur tour les ouverts $X_{\alpha}\times\mathbb{C}$
aux $\mathbb{A}^{1}$-fibrés triviaux $\rho_{\alpha}=p_{1}:X_{\alpha}\times\mathbb{C}\rightarrow X_{\alpha}$,
on conclut que $\rho:S\rightarrow X$ est un $\mathbb{A}^{1}$-fibré
devenant trivial sur le recouvrement ouvert canonique $\mathcal{U}=\left(X_{\alpha}\right)_{\alpha=\pm1,\,\pm2}$
de $X$. D'après (\ref{eq:partial_tilde_eqn}) ci-dessus, la restriction
à l'ouvert $S_{\alpha}\setminus C_{\alpha}\simeq\mathbb{C}^{*}\times\mathbb{C}$
de $S$ de la $\mathbb{C}_{+}$-action sur $S$ induite par la dérivation
$\partial$ s'identifie à celle induite par la $\mathbb{C}\left[x,x^{-1}\right]$
dérivation ${\displaystyle x^{2}\frac{\partial}{\partial f}.}$ de
$\mathbb{C}\left[x,x^{-1},y,z\right]$. Par conséquent, les isomorphismes
$\tau_{\alpha}:S_{\alpha}\stackrel{\sim}{\rightarrow}X_{\alpha}\times\mathbb{C}$
sont équivariants lorsque l'on munit les $X_{\alpha}\times\mathbb{C}$
des $\mathbb{C}_{+}$-actions par translations le long du second facteur.
Ainsi $\rho:S\rightarrow X$ est bien un fibré principal homogène
de groupe $\mathbb{C}_{+}$ et de base $X$, dont la classe d'isomorphisme
s'identifie à la classe de cohomologie $\left[g\right]\in\check{H}^{1}\left(\mathcal{U},\mathcal{O}_{X}\right)\simeq H^{1}\left(X,\mathcal{O}_{X}\right)$
du $1$-cocycle de \v{C}ech $\left\{ g_{\alpha\alpha'}\right\} \in C^{1}\left(\mathcal{U},\mathcal{O}_{X}\right)$. 

\end{enavant}

\begin{rem}
Une manière commode de coder les fonctions de transition $g_{\alpha\alpha'}\in\mathbb{C}\left[x,x^{-1}\right]$
définies ci-dessus est d'introduire l'arbre pondéré suivant : 

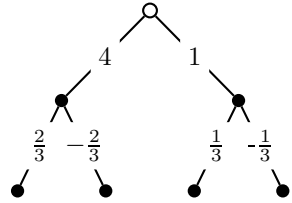
\begin{pspicture}(-3,1.6)(8,-1.5)

\rput(4,0){

\pstree[treemode=D,radius=2.5pt,treesep=1cm,levelsep=1.2cm]{\Tc{3pt}}{

\pstree{\TC*\mput*{$4$}}{\TC*\mput*{$\frac{2}{3}$}\TC*\mput*{$-\frac{2}{3}$}}

\pstree{\TC*\mput*{$1$}}{\TC*\mput*{$\frac{1}{3}$}\TC*\mput*{-$\frac{1}{3}$}}

}

}

\end{pspicture}

\noindent En effet, on peut associer aux branches de cet arbres les
polynômes $\sigma_{1}=1+x/3$, $\sigma_{-1}=1-x/3$, $\sigma_{2}=4+2x/3$
et $\sigma_{-2}=4-2x/3$ à partir desquels on retrouve les fonctions
de transitions ci-dessus en posant $g_{\alpha\alpha'}=x^{-2}\left(\sigma_{\alpha'}-\sigma_{\alpha}\right)\in\mathbb{C}\left[x,x^{-1}\right]$
pour tout $\alpha,\alpha'=\pm1,\pm2$. Cet exemple particulier est
en fait l'illustration d'une certaine correspondance plus générale
entre arbres pondérés et fibrés principaux homogènes au-dessus d'une
droite affine à plusieurs origines (cf. \cite{DubG03}). 
\end{rem}

\subsection{Comparaison de structures $\mathbb{A}^{1}$-fibrées}

\indent\newline\noindent Pour terminer la preuve du Théorème \ref{thm:Th_Principal},
il reste à montrer que la surface $S\subset\mathbb{C}^{3}$ d'équation
$xy-\left(f+1\right)\left(f+4\right)=0$, où $f=xz-y^{2}$, n'est
pas isomorphe à une surface de type $S_{Q,n}$ définie par une équation
de la forme $x^{n}z-Q\left(x,y\right)=0$. Admettons pour l'instant
le fait suivant qui sera établi ultérieurement :\\

\emph{Toute} $\mathbb{A}^{1}$-\emph{fibration} $\pi':S\rightarrow\mathbb{C}$
\emph{sur la surface} $S\subset\mathbb{C}^{3}$ \emph{définie par
l'équation} \[
x\left(xz^{2}-2y^{2}z+5z\right)+y^{4}-5y^{2}+4=0,\]
 \emph{est de la forme} $\pi'=\beta\circ\pi$, \emph{où} $\beta$
\emph{est un automorphisme de la droite affine} $\mathbb{C}$ \emph{et
où} $\pi={\rm pr}_{x}\mid_{S}:S\rightarrow\mathbb{C}$. 

\begin{enavant} Une surface $S_{Q,n}$ est elle aussi naturellement
munie d'une fibration ${\rm pr}_{x}:S_{Q,n}\rightarrow\mathbb{C}$
se restreignant en un $\mathbb{A}^{1}$-fibré trivial ${\rm pr}_{x}^{-1}\left(\mathbb{C}^{*}\right)\simeq\textrm{Spec}\left(\mathbb{C}\left[x,x^{-1}\right]\left[y\right]\right)\simeq\mathbb{C}^{*}\times\mathbb{C}$
au-dessus de $\mathbb{C}^{*}$. Cette fibration est surjective si
et seulement si $Q\left(0,y\right)$ est un polynôme non constant.
Si $Q\left(0,y\right)=0$, alors $S_{Q,n}$ est réductible, la fibre
${\rm pr}_{x}^{-1}\left(0\right)$ étant alors isomorphe à $\mathbb{C}^{2}$.
Sinon si $Q\left(0,y\right)=a\in\mathbb{C}^{*}$, alors la fibre ${\rm pr}_{x}^{-1}\left(0\right)$
est vide, de sorte que $S_{Q,n}$ est isomorphe à $\mathbb{C}^{*}\times\mathbb{C}$.
On en conclut que la surface $S$, qui est irréductible, ne peut être
isomorphe à une surface $S_{Q,n}$ du premier type. D'autre part,
il existe une obstruction de nature topologique à l'existence d'un
isomorphisme entre $S$ et une surface $S_{Q,n}$ du second type.
En effet, on sait que $S$ est isomorphe à une surface obtenue en
recollant quatre copies du plan affine en dehors d'un axe de coordonnées.
Il s'en suit de manière évidente que $H_{2}\left(S;\mathbb{Z}\right)\simeq\mathbb{Z}^{3}$,
de sorte que $S$ ne peut pas être isomorphe à $\mathbb{C}^{*}\times\mathbb{C}$.
Par conséquent, si $S$ est isomorphe à une surface du type $S_{Q,n}$
alors cette dernière est nécessairement définie par un polynôme $Q\left(x,y\right)$
tel que $Q\left(0,y\right)$ soit non constant. Dans ce cas, la fibration
induite ${\rm pr}_{x}:S_{Q,n}\rightarrow\mathbb{C}$ est une $\mathbb{A}^{1}$-fibration
dont l'unique fibre éventuellement non isomorphe à la droite affine
complexe est la fibre ${\rm pr}_{x}^{-1}\left(0\right)$. Puisque
nous avons supposé établi le fait que la fibration quotient $\pi={\rm pr}_{x}\mid_{S}:S\rightarrow\mathbb{C}$
est l'unique $\mathbb{A}^{1}$-fibration sur $S$ à automorphismes
de la base près, on conclut que s'il existe un isomorphisme $\phi:S\stackrel{\sim}{\rightarrow}S_{Q,n}$
alors il existe un automorphisme $\zeta$ de $\mathbb{C}$ fixant
l'origine tel $\pi\circ\zeta={\rm pr}_{x}\mid_{S_{Q,n}}\circ\phi$.
La proposition suivante termine par conséquent la preuve du Théorème
\ref{thm:Th_Principal}. 

\end{enavant} 

\begin{prop}
Il n'existe aucun isomorphisme $\phi$ de ce type. 
\end{prop}
\begin{proof}
Pour qu'un isomorphisme $\phi$ comme ci-dessus puisse exister, il
est nécessaire que la fibre ${\rm pr}_{x}^{-1}\left(0\right)$ de
la $\mathbb{A}^{1}$-fibration ${\rm pr}_{x}:S_{Q,n}\rightarrow\mathbb{C}$
soit réduite, constituée d'exactement quatre composantes irréductibles,
toutes isomorphes à la droite affine complexe. Cela impose donc que
$Q\left(0,y\right)$ soit un polynôme de degré quatre, à racines simples.
Dans ces conditions, la fonction $y$ se restreint sur $S_{Q,h}$
en une fonction induisant une coordonnée sur les fibres de ${\rm pr}_{x}:S_{Q,n}\rightarrow\mathbb{C}$
au-dessus de $\mathbb{C}^{*}$ et localement constante sur le fibre
${\rm pr}_{x}^{-1}\left(0\right)$. De plus $y\mid_{S_{Q,n}}$ distingue
les composantes irréductibles de ${\rm pr}_{x}^{-1}\left(0\right)$.
S'il existait un isomorphisme $\phi:S\stackrel{\sim}{\rightarrow}S_{Q,n}$
alors $y_{\phi}=y\mid_{S_{Q,n}}\circ\phi$ serait une fonction sur
$S$ jouissant des mêmes propriétés. Dans les trivialisations locales
$\tau_{\alpha}^{*}:\mathbb{C}\left[x\right]\left[u_{\alpha}\right]\stackrel{\sim}{\rightarrow}\mathcal{O}\left(S_{\alpha}\right)$
construites en \ref{txt:Trivilisations_locales} ci-dessus, $\tau_{\alpha}^{*}y_{\phi}$
devrait donc s'écrire sous la forme \[
\left(\tau_{\alpha}^{*}\right)^{-1}y_{\phi}=x^{n_{\alpha}}u_{\alpha}+P_{\alpha}\left(x\right)\in\mathbb{C}\left[x\right]\left[u_{\alpha}\right]\]
 où $n_{\alpha}\geq1$ pour tout $\alpha=\pm1,\pm2$ et où $P_{\alpha}\left(0\right)\neq P_{\alpha'}\left(0\right)$
pour tout $\alpha\neq\alpha'$. Puisque la fonction $y_{\phi}$ est
définie globalement sur $S$, on devrait donc avoir \[
x^{n_{\alpha'}}u_{\alpha'}+P_{\alpha'}\left(x\right)=\left(\tau_{\alpha'}^{*}\right)^{-1}\tau_{\alpha}^{*}\left(x^{n_{\alpha}}u_{\alpha}+P_{\alpha}\left(x\right)\right)=x^{n_{\alpha}}\left(u_{\alpha'}+g_{\alpha\alpha'}\right)+P_{\alpha}\left(x\right)\]
 dans $\mathbb{C}\left[x\right]\left[u_{\alpha'}\right]$ pour tout
$\alpha'\neq\alpha$, où $g_{\alpha\alpha'}\in\mathbb{C}\left[x,x^{-1}\right]$
désigne la fonction de transition construite en \ref{txt:Fonctions_de_transition}
ci-dessus. On en déduit que $n_{\alpha}=n_{\alpha'}=n\geq2$ pour
tout $\alpha,\alpha'\in\left\{ \pm1,\,\pm2\right\} $. Si $n>2$ alors
$x^{n}\left(u_{\alpha'}+g_{\alpha\alpha'}\right)\in x\mathbb{C}\left[x\right]\left[u_{\alpha'}\right]$
de sorte que l'égalité ci-dessus conduirait à $P_{\alpha}\left(0\right)=P_{\alpha'}\left(0\right)$
pour tout $\alpha\neq\alpha'$ en contradiction avec notre hypothèse.
Ainsi $n=2$, et l'on tire de l'égalité ci-dessus que $P_{\alpha'}\left(0\right)=\left(\left(\alpha'\right)^{2}-\alpha^{2}\right)+P_{\alpha}\left(0\right)$
pour tout $\alpha\neq\alpha'$. Cela conduit de nouveau à une contradiction
puisque pour $\alpha'=-\alpha$, on trouve que $P_{\alpha'}\left(0\right)=P_{\alpha}\left(0\right)$.
Il ne peut donc pas exister d'isomorphisme $\phi:S\stackrel{\sim}{\rightarrow}S_{Q,h}$
comme ci-dessus. 
\end{proof}

\subsection{Unicité à isomorphisme près de la fibration $\pi:S\rightarrow S/\!/\mathbb{C}_{+}$ }

\indent\newline\noindent Le but de ce paragraphe est d'établir le
résultat admis auparavant, à savoir le fait que la fibration quotient
$\pi:S\rightarrow S/\!/\mathbb{C}_{+}\simeq\mathbb{C}$ est l'unique
$\mathbb{A}^{1}$-fibration sur $S$ à automorphismes de la base près.
En se basant sur des résultats antérieurs de M.H. Gizatullin \cite{Giz71},
J. Bertin \cite{Ber83} à montré que cela était le cas dès lors que
$S$ admet une compactification projective minimale $\bar{S}$ pour
laquelle le diviseur de bord $\bar{S}\setminus S$ n'est pas une chaîne
de courbes rationnelles propres non singulières. Cette caractérisation
repose essentiellement sur l'idée simple que la structure du bord
d'une compactification d'une surface affine réglée $S$ est intimement
liée à la nature des singularités à l'infini des pinceaux de courbes
rationnelles sur $S$. Malheureusement, la technicité des preuves
données dans \emph{loc. cit.,} ainsi d'ailleurs que dans d'autres
articles exposant des résultats du même type (cf. \cite{DaRus02},
\cite{Dub02} et \cite{GuMiy05}), rend sa lecture difficile. C'est
pourquoi nous donnons ici, en guise d'illustration de ces méthodes
générales, une preuve complète de ce résultat dans le cas particulier
de la surface $S$ du Théorème \ref{thm:Th_Principal}. Plus précisément,
nous allons établir le fait suivant : 

\begin{prop}
\label{pro:Unicit=E9_des_fibrations} Toute $\mathbb{A}^{1}$-fibration
$\pi':S\rightarrow\mathbb{C}$ sur la surface $S\subset\mathbb{C}^{3}$
définie par l'équation \[
x\left(xz^{2}-2y^{2}z+5z\right)+y^{4}-5y^{2}+4=0,\]
 est de la forme $\pi'=\beta\circ\pi$, où $\beta$ est un automorphisme
de la droite affine $\mathbb{C}$ et où $\pi={\rm pr}_{x}\mid_{S}:S\rightarrow\mathbb{C}$
désigne la $\mathbb{A}^{1}$-fibration quotient associée à la $\mathbb{C}_{+}$-action
sur $S$ induite par la dérivation $\partial$ du Théorème \ref{thm:Th_Principal}. 
\end{prop}
\begin{enavant} La preuve de la Proposition \ref{pro:Unicit=E9_des_fibrations}
se décompose en deux parties. Nous construisons tout d'abord, par
éclatements successifs au-dessus $\mathbb{P}^{1}\times\mathbb{P}^{1}$,
une surface projective lisse $\bar{S}$ munie d'une $\mathbb{P}^{1}$-fibration
$\bar{\pi}:\bar{S}\rightarrow\mathbb{P}^{1}$ compactifiant $\pi:S\rightarrow\mathbb{C}$
au sens qu'il existe un diagramme commutatif \[\xymatrix{ S \ar[d]_{\pi} \ar[r]^{i_S} & \bar{S} \ar[d]_{\bar{\pi}} \\ \mathbb{C} \ar[r]^{i} & \mathbb{P}^1}\]

\noindent dans lequel $i_{S}:S\hookrightarrow\bar{S}$ et $i:S/\!/\mathbb{C}_{+}\hookrightarrow\mathbb{P}^{1}$
sont des immersions ouvertes. Nous montrons ensuite que la structure
du diviseur de bord $B=\bar{S}\setminus S$ interdit l'existence de
$\mathbb{A}^{1}$-fibrations $p:S\rightarrow\mathbb{C}$ sur $S$
dont les fibres générales seraient transverses à celles de $\pi:S\rightarrow\mathbb{C}$. 

\end{enavant}

\subsubsection{Compactification de la fibration quotient}

\indent\newline\noindent Toute surface affine irréductible $S$ munie
d'une $\mathbb{A}^{1}$-fibration $p:S\rightarrow\mathbb{C}$ admet
des morphismes birationnels $\left(p,\psi\right):S\rightarrow\mathbb{C}\times\mathbb{C}$.
En effet, étant donné un ouvert affine $U\subset\mathbb{C}$ au-dessus
duquel $p$ devient un $\mathbb{A}^{1}$-fibré trivial, toute fonction
régulière $\psi\in\mathcal{O}\left(S\right)$ induisant un isomorphisme
$p^{-1}\left(U\right)\simeq\textrm{Spec}\left(\mathcal{O}\left(U\right)\left[\psi\right]\right)$
permet de définir un morphisme birationnel du type annoncé. La structure
des morphismes birationnels $\sigma:V'\rightarrow V$ entre variétés
affines, souvent appelés \emph{modifications affines} \cite{KaZa99},
est connue. Tout morphisme de ce type est obtenu en faisant d'abord
éclater un sous-schéma fermé $Y$ de $V$, puis en rejetant la transformée
propre d'un diviseur de $V$ contenant $Y$. Lorsque les fibres de
la $\mathbb{A}^{1}$-fibration $p:S\rightarrow\mathbb{C}$ sont réduites,
un procédé introduit par K.H. Fieseler \cite{Fie94} et précisé par
l'auteur \cite{DubG03} permet d'obtenir une factorisation des morphismes
birationnels ci-dessus en une succession explicite d'éclatements ponctuels.
C'est ce procédé que nous employons ci-dessous pour construire une
compactification $\bar{S}$ de la surface $S=\textrm{Spec}\left(\mathbb{C}\left[x,y,z\right]/\left(xy-\left(f+1\right)\left(f+4\right)\right)\right)$
du Théorème \ref{thm:Th_Principal}.

\begin{enavant} Pour débuter la construction, il nous faut choisir
un morphisme birationnel $\sigma=\left(\pi,\psi\right):S\rightarrow\mathbb{C}\times\mathbb{C}$.
Or d'après la Proposition \ref{pro:Fibration_sur_S}, la restriction
du polynôme $f=xz-y^{2}$ sur $S$ induit un isomorphisme $\pi^{-1}\left(\mathbb{C}^{*}\right)\simeq\textrm{Spec}\left(\mathbb{C}\left[x,x^{-1}\right]\left[f\right]\right)$.
Ainsi, en prenant $\psi=f$, on obtient un morphisme birationnel $\sigma=\left(\pi,f\right):S\rightarrow\mathbb{C}\times\mathbb{C}$,
$\left(x,y,z\right)\mapsto\left(x,xz-y^{2}\right)$ induisant un isomorphisme
$\pi^{-1}\left(\mathbb{C}^{*}\right)\stackrel{\sim}{\rightarrow}\mathbb{C}^{*}\times\mathbb{C}\subset\mathbb{C}\times\mathbb{C}$.
D'autre part, il est clair que l'image de $S$ par le morphisme $\phi:\mathbb{C}^{3}\rightarrow\mathbb{C}^{3}$,
$\left(x,y,z\right)\mapsto\left(x,xz-y^{2},y\right)$ est contenue
dans la surface $S'\subset\mathbb{C}^{3}=\textrm{Spec}\left(\mathbb{C}\left[x,y',z'\right]\right)$
d'équation $xz'-\left(y'+1\right)\left(y'+4\right)=0$. De plus le
morphisme birationnel $\sigma$ se factorise via $\phi$ par la projection
$\sigma'={\rm pr}_{x,y'}\mid_{S'}:S'\rightarrow\mathbb{C}\times\mathbb{C}$,
$\left(x,y',z'\right)\mapsto\left(x,y'\right)$. Cette surface $S'$
est munie d'une $\mathbb{A}^{1}$-fibration $\pi'={\rm pr}_{x}\mid_{S'}:S'\rightarrow\mathbb{C}$,
se restreignant en un $\mathbb{A}^{1}$-fibré trivial $\mathbb{C}^{*}\times\mathbb{C}$
au-dessus de $\mathbb{C}^{*}$, avec $y'\mid_{S'}$ comme coordonnée
naturelle sur le second facteur. La fibre $\left(\pi'\right)^{-1}\left(0\right)$
est constituée de la réunion disjointe de deux droites affines complexes
$C'_{-\alpha^{2}}$, $\alpha=\pm1,\,\pm2$ d'équations respectives
$\left\{ x=0,\, y'=-1\right\} $ et $\left\{ x=0,\, y'=-4\right\} $
dont les images par la projection birationnelle $\sigma':S'\rightarrow\mathbb{C}\times\mathbb{C}$
sont les deux points $p_{-1}=\left(0,-1\right)$ et $p_{-4}=\left(0,-4\right)$
respectivement. Par construction, l'image par $\phi$ d'une composante
$C_{\alpha}$, $\alpha=\pm1,\,\pm2$, de $\pi^{-1}\left(0\right)$
est constituée du point $q_{\alpha}=\left(0,-\alpha^{2},\alpha\right)\in C'_{-\alpha^{2}}$. 

\end{enavant}

\begin{enavant} Il se trouve que la surface $S'\subset\mathbb{C}^{3}$
admet tout comme $S$ une structure de $\mathbb{A}^{1}$-fibré $\rho':S'\rightarrow X'$,
cette fois au-dessus de la droite affine à deux origines. En effet,
les ouverts $S'_{-1}=S'\setminus C'_{-4}$ et $S'_{-4}=S'\setminus C'_{-1}$
sont tous deux isomorphes à $\mathbb{C}\times\mathbb{C}$ muni respectivement
des ``coordonnées''\[
\left(x,v_{-\alpha^{2}}={\displaystyle \frac{y'+\alpha^{2}}{x}\mid_{\tilde{S}}}\right)\in\mathcal{O}\left(S'_{-\alpha^{2}}\right)^{2},\quad\alpha^{2}=-1,\,-4.\]
 La restriction sur $S'_{-\alpha^{2}}$ de la projection birationnelle
$\sigma':S'\rightarrow\mathbb{C}\times\mathbb{C}$ coïncide alors
avec le morphisme $\sigma'_{-\alpha^{2}}:S'_{-\alpha^{2}}\rightarrow\mathbb{C}^{2}$,
$\left(x,v_{-\alpha^{2}}\right)\mapsto\left(x,xv_{-\alpha^{2}}-\alpha^{2}\right)$,
qui n'est autre que l'expression dans l'une des deux cartes affines
usuelles de l'éclatement du point $p_{-\alpha^{2}}$ de $\mathbb{C}^{2}$.
Par construction, la préimage par $\phi\mid_{S}^{S'}$ de l'ouvert
$S'_{-\alpha^{2}}$ de $S'$ est constituée de la réunion des ouverts
$S_{\alpha}$ et $S_{-\alpha}$ de $S$. De plus, pour tout $\alpha=\pm1,\,\pm2$,
le morphisme induit $\phi\mid_{S}^{S'}:S_{\alpha}\rightarrow S'_{-\alpha^{2}}$
coïncide de nouveau dans les coordonnées $\left(x,u_{\alpha}\right)$
et $\left(x,v_{-\alpha^{2}}\right)$ introduites ci-dessus avec l'éclatement
$\left(x,u_{\alpha}\right)\mapsto\left(x,xu_{\alpha}-{\displaystyle \left(-1\right)^{\left|\alpha\right|}\frac{\alpha}{3}}\right)$
du point $q_{\alpha}$ de $S'_{-\alpha^{2}}\simeq\mathbb{C}^{2}$.
La situation peut donc se résumer par la figure suivante : \psset{unit=0.4cm}

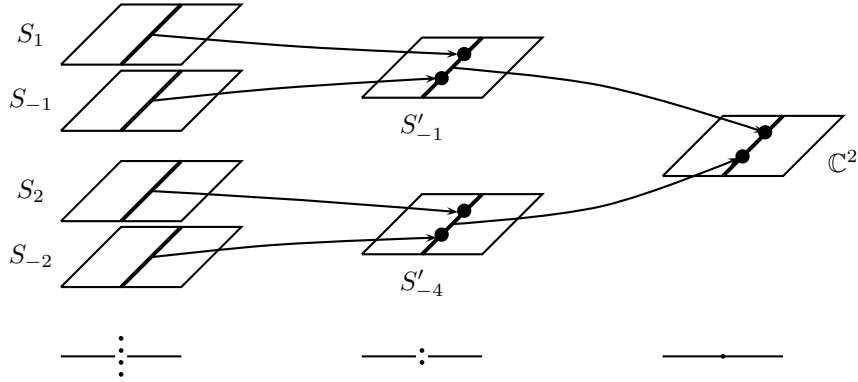
\newsavebox{\para}  

\savebox{\para}(5,4){\psline(0,0)(2,2)(6,2)(4,0)(0,0) \psline[linewidth=1.5pt] (2,0)(4,2)}

\newsavebox{\parapoint}

\savebox{\parapoint}(5,4){\rput(0,0){\usebox{\para}}

\psdots[dotscale=1.5](3.4,1.45)(2.65,0.65)}

\begin{pspicture}(-1.5,7)(10,-7)

\rput(4,3.7){\usebox{\para}}

\rput(3,4.7){$S_{1}$}

\rput(4,1.5){\usebox{\para}}

\rput(3,2.5){$S_{-1}$}

\rput(4,-1.5){\usebox{\para}}

\rput(3,-0.5){$S_{2}$}

\rput(4,-3.7){\usebox{\para}}

\rput(3,-2.7){$S_{-2}$}

\rput(14,2.6){\usebox{\parapoint}}

\rput(16,1.6){$S'_{-1}$}

\rput(14,-2.6){\usebox{\parapoint}}

\rput(16,-3.6){$S'_{-4}$}

\rput(24,0){\usebox{\parapoint}}

\rput(30,0.5){$\mathbb{C}^2$}

\pscurve{->}(7,4.7)(12,4.3)(17.2,4.05)

\pscurve{->}(7,2.5)(11,2.9)(16.5,3.25)

\pscurve{->}(7,-0.5)(12,-0.8)(17.2,-1.2)

\pscurve{->}(7,-2.7)(11,-2.3)(16.5,-2.05)

\pscurve{->}(17,3.6)(22,3)(27.4,1.45)

\pscurve{->}(17,-1.6)(22,-1)(26.55,0.6)

\psline(4,-6)(5.8,-6)\psline(6.2,-6)(8,-6)

\psdots[dotscale=0.5](6,-5.4)(6,-5.8)(6,-6.2)(6,-6.6)

\psline(14,-6)(15.8,-6)\psline(16.2,-6)(18,-6)

\psdots[dotscale=0.5](16,-5.8)(16,-6.2)

\psline(24,-6)(28,-6)\psdots[dotscale=0.5](26,-6)

\end{pspicture}

\noindent dans laquelle chaque flèche réprésente l'unique endomorphisme
affine birationnel de $\mathbb{C}^{2}$ contractant la courbe à la
source sur le point au but et induisant un isomorphisme entre $\mathbb{C}^{2}$
privé de cette courbe et son image. 

\end{enavant}

\begin{enavant} Ainsi, à partir de toute compactification du plan
affine $\mathbb{C}^{2}$ en une surface projective lisse $\bar{V}$,
on obtiendra une compactification de $S$ en une surface projective
lisse obtenue en éclatant tout d'abord les points $p_{-1}$ et $p_{-4}$
puis les points infiniment voisins $q_{\pm1}$ et $q_{\pm2}$ qui
se situent sur les diviseurs exceptionnels au-dessus de $p_{-1}$
et $p_{-4}$ respectivement. En particulier, nous pouvons considérer
une compactification de $\mathbb{C}^{2}=\mathbb{C}\times\mathbb{C}$
sous la forme d'une surface de Hirzebruch $\bar{\pi}_{0}:\mathbb{F}_{n}=\mathbb{P}\left(\mathcal{O}_{\mathbb{P}^{1}}\oplus\mathcal{O}_{\mathbb{P}^{1}}\left(-n\right)\right)\rightarrow\mathbb{P}^{1}$
de telle sorte que $\bar{\pi}_{0}$ prolonge la première projection
${\rm pr}_{1}:\mathbb{C}\times\mathbb{C}\rightarrow\mathbb{C}$. Si
$n\geq1$, nous supposerons de plus que $\mathbb{F}_{n}\setminus\mathbb{C}^{2}$
est la réunion de la section négative canonique $D$ de $\bar{\pi}_{0}$
et d'une fibre $F_{\infty}$ de $\bar{\pi}_{0}$ au dessus de $\infty=\mathbb{P}^{1}\setminus\mathbb{C}$.
Sinon, si $n=0$, alors $\mathbb{F}_{0}\setminus\mathbb{C}^{2}=D\cup F_{\infty}$
est constitué de la réunion d'une fibre de chacun des réglages évidents
sur $\mathbb{F}_{0}\simeq\mathbb{P}^{1}\times\mathbb{P}^{1}$. Nous
notons $F_{0}\simeq\mathbb{P}^{1}$ la fibre de $\bar{\pi}_{0}$ au
dessus de $0\in\mathbb{C}\subset\mathbb{P}^{1}$, $E_{-1}$ et $E_{-4}$
les diviseurs exceptionnels de l'éclatement $\bar{\sigma}':\bar{S}'\rightarrow\mathbb{F}_{n}$
des points $p_{-1}$ et $p_{-4}$ de $\mathbb{C}^{2}\subset\mathbb{F}_{n}$.
Nous désignons par $\beta:\bar{S}\rightarrow\bar{S}'$ l'éclatement
des points $q_{\alpha}$, $\alpha=\pm1,\,\pm2$ de $S'\subset\bar{S}'$,
les diviseurs exceptionnels correspondant étant les adhérences $\bar{C}_{\alpha}$
dans $\bar{S}$ des composantes $C_{\alpha}$, $\alpha=\pm1,\,\pm2$,
de la fibre $\pi^{-1}\left(0\right)$. Par construction, nous obtenons
une compactification $\bar{S}$ de $S$ munie d'une $\mathbb{P}^{1}$-fibration
$\bar{\pi}:\bar{S}\rightarrow\mathbb{P}^{1}$ relevant $\bar{\pi}_{0}:\mathbb{F}_{n}\rightarrow\mathbb{P}^{1}$
et prolongeant la $\mathbb{A}^{1}$-fibration $\pi:S\rightarrow\mathbb{C}$.
Le bord $B=\bar{S}\setminus S$ est un diviseur à croisements normaux
simples, dont les composantes irréductibles sont les transformées
strictes dans $\bar{S}$ des courbes $D$, $F_{\infty}$, $F_{0}$,
$E_{-1}$ et $E_{-4}$, que nous noterons encore par les mêmes symboles.
Avec ces conventions, les indices d'auto-intersection des composantes
de $B$ sont les suivants : $\left(D^{2}\right)=-n$, $\left(F_{\infty}^{2}\right)=0$,
$\left(F_{0}^{2}\right)=-2$, $\left(E_{-1}^{2}\right)=\left(E_{-4}^{2}\right)=-3$.
Les adhérences $\bar{C}_{\alpha}$ des composantes de $\pi^{-1}\left(0\right)$
sont quant à elles des courbes exceptionnelles de première espèce.
Cette compactification $\bar{S}$ de $S$ est minimale, au sens qu'aucune
des composantes irréductibles du bord $B$ n'est une courbe exceptionnelle
de première espèce. Le graphe dual du diviseur $B\cup\bar{C}_{1}\cup\bar{C}_{-1}\cup\bar{C}_{2}\cup\bar{C}_{-2}$
est le suivant :

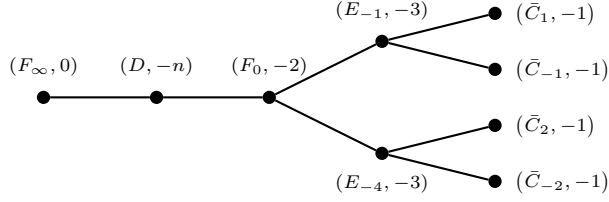
\begin{pspicture}(-3,1.8)(8,-1.5)

{\scriptsize

\rput(4,0){

\pstree[treemode=R,radius=2.5pt,treesep=0.45cm,levelsep=1.5cm]{\Tc*{2.5pt} ~[tnpos=a]{$\left(F_{\infty},0\right)$}} {

    \pstree{\TC*~[tnpos=a]{$\left(D,-n\right)$}}{

        \pstree{\TC*~[tnpos=a]{$\left(F_0,-2\right)$}}{

           \pstree{\TC*~[tnpos=a]{$\left(E_{-1},-3\right)$}}{\TC*~[tnpos=r]{$\left(\bar{C}_1,-1\right)$} \TC*~[tnpos=r]{$\left(\bar{C}_{-1},-1\right)$}}

           \pstree{\TC*~[tnpos=b]{$\left(E_{-4},-3\right)$}}{\TC*~[tnpos=r]{$\left(\bar{C}_2,-1\right)$}\TC*~[tnpos=r]{$\left(\bar{C}_{-2},-1\right)$}}

}

}

}

}

}

\end{pspicture}

\end{enavant}

\subsubsection{Non existence de $\mathbb{A}^{1}$-fibrations transverses}

\indent\newline\noindent D'après ce qui précède, le graphe dual du
diviseur de bord $B=\bar{S}\setminus S$ de la compactification minimale
$\bar{S}$ de $S$ construite ci-dessus n'est pas une chaîne. Suivant
un critère dû à J. Bertin \cite{Ber83}, l'existence d'une compactification
possédant cette propriété suffit à garantir que la fibration naturelle
$\pi:S\rightarrow\mathbb{C}$ est l'unique $\mathbb{A}^{1}$-fibration
sur $S$ à automorphismes de la base près. Bien qu'il puisse sembler
a priori mystérieux, ce critère ne repose au final que sur le fait
bien connu que toute surface projective lisse $\bar{V}$ munie d'une
$\mathbb{P}^{1}$-fibration $p:\bar{V}\rightarrow\mathbb{P}^{1}$
est obtenue en effectuant une suite d'éclatements de centres lisses
à partir d'une surface de Hirzebruch $\mathbb{F}_{n}\rightarrow\mathbb{P}^{1}$.
Afin de donner au lecteur un aperçu de la méthode utilisée dans la
preuve du critère de Bertin, nous allons établir directement qu'à
automorphismes de la base près, $S$ n'admet pas d'autre $\mathbb{A}^{1}$-fibration
que $\pi:S\rightarrow\mathbb{C}$. 

\begin{enavant} Rappelons le résultat classique de théorie de surfaces
(cf. Chapitre 4 de \cite{GrHar78} ) selon lequel toute surface projective
normale $\bar{V}$ munie d'une $\mathbb{P}^{1}$-fibration $p:\bar{V}\rightarrow\mathbb{P}^{1}$
est obtenue par une suite d'éclatements de centres lisses à partir
d'une surface de Hirzebruch $\mathbb{F}_{n}\rightarrow\mathbb{P}^{1}$
convenable, de telle manière que $p$ relève la structure naturelle
de $\mathbb{P}^{1}$-fibré sur $\mathbb{F}_{n}$. En particulier,
toute fibre d'une $\mathbb{P}^{1}$-fibration sur une surface projective
lisse est supportée par un diviseur à croisements normaux simples
dont les composantes irréductibles sont des courbes rationnelles propres
et lisses. De plus, les indices d'auto-intersection des composantes
irréductibles d'une fibre non irréductible sont tous négatifs, l'une
au moins des ces composantes étant nécessairement une courbe exceptionnelle
de première espèce. 

\end{enavant}

\begin{enavant} Supposons que $S$ soit munie d'une seconde $\mathbb{A}^{1}$-fibration
$\pi':S\rightarrow\mathbb{C}$ dont les fibres générales sont distinctes
de celles de $\pi$. En particulier, $\pi$ (resp. $\pi'$) n'est
pas constant sur les fibres générales de $\pi'$ (resp. de $\pi$).
Nous pouvons naturellement considérer $\pi'$ comme une application
rationnelle $\bar{\pi}':\bar{S}\dashrightarrow\mathbb{P}^{1}=\mathbb{C}\cup\left\{ \infty'\right\} $
sur une des compactifications $\bar{S}$ de $S$ construite ci-dessus.
Supposons pour fixer les idées que $\bar{S}=\bar{S}_{0}$ est la compactification
obtenue à partir de $\mathbb{F}_{0}=\mathbb{P}^{1}\times\mathbb{P}^{1}$
par la suite d'éclatements décrite précédemment. On peut interpréter
$\pi'$ comme un pinceau de courbes rationnelles sur $\bar{S}$ dont
les membres généraux sont les adhérences $\bar{F}'\simeq\mathbb{P}^{1}$
dans $\bar{S}$ des fibres générales $F'\simeq\mathbb{A}^{1}$ de
$\pi'$. Le point $\bar{F}'\setminus F'$ est nécessairement contenu
dans la fibre $F_{\infty}=\bar{\pi}^{-1}\left(\infty\right)$ car
sinon $\bar{\pi}$ serait bornée, donc constante, sur les fibres générales
de $\pi'$, ce qui contredirait notre hypothèse. On en déduit en particulier
que si le pinceau $\bar{\pi}'$ a un point de base $p_{0}$ (nécessairement
unique), alors ce dernier est contenu dans $F_{\infty}$. L'application
$\bar{\pi}'$ est définie sur la section $D$ de $\bar{\pi}$, privée
éventuellement du point $p_{0}$, et $\bar{\pi}\left(D\setminus\left\{ p_{0}\right\} \right)=\infty'$
car sinon $\bar{\pi}$ serait bornée, donc constante, sur les fibres
générales de $\pi$, ce qui contredirait de nouveau l'hypothèse. Puisque
$F_{0}\cup E_{-1}\cup E_{-4}$ est disjoint de $F_{\infty}$ et ne
contient pas l'éventuel point base $p_{0}$ de $\bar{\pi}'$, on déduit
que la restriction de $\bar{\pi}'$ sur $G=F_{0}\cup E_{-1}\cup E_{-4}\cup\left(D\setminus\left\{ p_{0}\right\} \right)$
est constante, de valeur $\infty'$. Par conséquent, après résolution
de l'éventuel point base $p_{0}$ par une succession d'éclatements
de centres lisses $\sigma:\tilde{S}\rightarrow\bar{S}$, la transformée
propre de $\bar{G}'$ est contenue dans la fibre $\tilde{\pi}^{-1}\left(\infty'\right)$
de la $\mathbb{P}^{1}$-fibration sur $\tilde{S}$ relevant $\bar{\pi}':\bar{S}\rightarrow\mathbb{P}^{1}$.
La figure suivante décrit l'allure du graphe dual de $\bar{G}'$ ainsi
que les indices d'auto-intersection de chacune de ses composante dans
les trois cas possibles, à savoir $\bar{\pi}'$ n'a pas de point d'indétermination,
l'unique point d'indétermination $p_{0}$ de $\bar{\pi}'$ appartient
à $F_{\infty}\setminus D$ , ou $p_{0}=F_{\infty}\cap D$: 

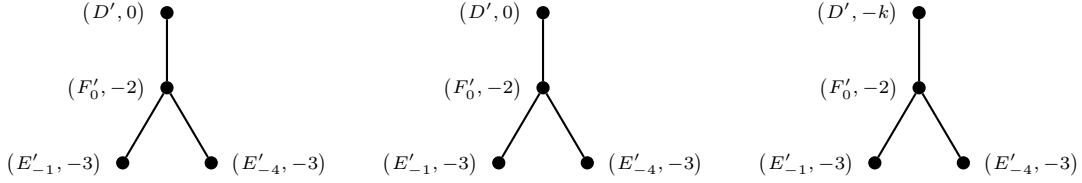
\begin{pspicture}(2,1.8)(8,-1.5)

{\scriptsize

\rput(4,0){

\pstree[treemode=D,radius=2.5pt,treesep=1cm,levelsep=1cm]{\Tc*{2.5pt}~[tnpos=l]{$\left(D',0\right)$} }{

 \pstree{\TC*~[tnpos=l]{$\left(F_0',-2\right)$}}{\TC*~[tnpos=l]{$\left(E_{-1}',-3\right)$}\TC*~[tnpos=r]{$\left(E_{-4}',-3\right)$}}

}

}

\rput(9,0){

\pstree[treemode=D,radius=2.5pt,treesep=1cm,levelsep=1cm]{\Tc*{2.5pt}~[tnpos=l]{$\left(D',0\right)$} }{

 \pstree{\TC*~[tnpos=l]{$\left(F_0',-2\right)$}}{\TC*~[tnpos=l]{$\left(E_{-1}',-3\right)$}\TC*~[tnpos=r]{$\left(E_{-4}',-3\right)$}}

}

} 

\rput(14,0){

\pstree[treemode=D,radius=2.5pt,treesep=1cm,levelsep=1cm]{\Tc*{2.5pt}~[tnpos=l]{$\left(D',-k\right)$} }{

 \pstree{\TC*~[tnpos=l]{$\left(F_0',-2\right)$}}{\TC*~[tnpos=l]{$\left(E_{-1}',-3\right)$}\TC*~[tnpos=r]{$\left(E_{-4}',-3\right)$}}

}

} 

}

\end{pspicture}

\noindent Dans les deux premiers cas, $\tilde{\pi}^{-1}\left(\infty'\right)$
serait réductible tout en contenant une composante d'indice d'auto-intersection
$0$, à savoir la transformée stricte de $D$, ce qui est impossible.
Dans le dernier cas, l'indice d'auto-intersection de la transformée
stricte de $D$ est au plus égal à $-1$. La fibre $\tilde{\pi}^{-1}\left(\infty'\right)$
devant pouvoir être contractée sur une courbe rationnelle lisse par
une succession de contractions de $\left(-1\right)$-courbes, il sera
nécessaire de contracter à une certaine étape de l'opération au moins
$3$ courbes parmis les transformées strictes $D'$, $F_{0}'$, $E_{-1}'$
et $E_{-4}'$ de $D$, $F_{0}$, $E_{-1}$ et $E_{-4}$ respectivement.
Ces trois dernières ayant un indice d'auto-intersection inférieur
ou égal à $-2$, toute contraction de ce type devra débuter par la
contraction de $D'$ suivie de la contraction de l'image de $F'_{0}$
qui sera devenue une $\left(-1\right)$-courbe. L'image de $\tilde{\pi}^{-1}\left(\infty'\right)$
par cette succession de contractions ne peut pas se réduire aux images
de $E'_{-1}$ et $E_{-4}'$ par cette seconde contraction. En effet,
ces dernières seraient alors des $\left(-2\right)$-courbes transverses,
qui, d'après ce qui précède, ne peuvent constituer à elles seules
une fibre d'une $\mathbb{P}^{1}$-fibration . D'autre part, si l'image
de $\tilde{\pi}^{-1}\left(\infty'\right)$ contient une autre composante
que les images de $E_{-1}'$ et $E_{-4}'$, alors elle contient en
particulier une composante qui intersecte ces courbes en leur unique
point d'intersection. Cela est de nouveau impossible puisqu'une fibre
d'une $\mathbb{P}^{1}$-fibration est supportée par un diviseur à
croisements normaux simples. 

\end{enavant}

\begin{enavant} La discussion ci-dessus montre que les fibres générales
de toute $\mathbb{A}^{1}$-fibration $\pi':S\rightarrow\mathbb{C}$
sur la surface $S$ du Théorème \ref{thm:Th_Principal} doivent coïncider
avec les fibres générales de la fibration naturelle $\pi={\rm pr}_{x}:S\rightarrow\mathbb{C}$,
donc, plus généralement, que toute fibre de $\pi'$ est contenue dans
une fibre de $\pi$. Par conséquent, toute $\mathbb{A}^{1}$-fibration
$\pi'$ sur $S$ se prolonge en une $\mathbb{P}^{1}$-fibration $\bar{\pi}':\bar{S}\rightarrow\mathbb{P}^{1}=\mathbb{C}\cup\left\{ \infty'\right\} $
sur les compactifications $\bar{S}$ de $S$ construites au paragraphe
ci-dessus, dont les fibres coïncident avec celle de $\bar{\pi}:\bar{S}\rightarrow\mathbb{P}^{1}$.
La courbe $D\simeq\mathbb{P}^{1}$ s'identifie alors à une section
commune de $\bar{\pi}'$ de $\bar{\pi}$, ce qui implique qu'il existe
un isomorphisme $\phi:\mathbb{P}^{1}\stackrel{\sim}{\rightarrow}\mathbb{P}^{1}$
tel que $\bar{\pi}'=\phi\circ\bar{\pi}$. Puisque l'on a nécessairement
$\phi\left(\infty\right)=\infty'$, $\phi$ induit un isomorphisme
$\beta:\mathbb{C}\stackrel{\sim}{\rightarrow}\mathbb{C}$ tel que
$\pi'=\beta\circ\pi$, ce qui conclut la preuve de la Proposition
\ref{pro:Unicit=E9_des_fibrations}. 

\end{enavant}

\section{Application aux automorphismes de $\mathbb{C}^{3}$ }

L'existence de la surface $S$ considérée, et surtout de la $\mathbb{C}_{+}$-action
sur $\mathbb{C}^{3}$ qui l'accompagne, a des répercussions sur la
nature de groupe d'automorphismes de $\mathbb{C}^{3}$. Nous allons
en particulier montrer comment construire à partir de la dérivation
localement nilpotente $\partial$ du Théorème \ref{thm:Th_Principal}
un automorphisme sauvage de l'espace affine dont, à la connaissance
de l'auteur, il n'existe aucune mention antérieure dans la littérature.

\subsection{Généralités sur les automorphismes algébriques de $\mathbb{C}^{3}$}

\indent\newline\noindent Il est évident que tout automorphisme algébrique
de la droite affine complexe est une transformation affine de la forme
$x\mapsto ax+b$, où $a\in\mathbb{C}^{*}$ et ou $b\in\mathbb{C}$.
La structure du groupe d'automorphismes ${\rm Aut}\left(\mathbb{C}^{2}\right)$
du plan affine est déjà notoirement plus complexe. Outre les automorphismes
affines généralisant ceux du type précédent, de nouveaux automorphismes,
dits \emph{triangulaires}, font leur apparition. Rappelons qu'un automorphisme
algébrique $\phi$ de $\mathbb{C}^{2}$ est dit triangulaire (ou plutôt
triangulable) s'il existe un système de coordonnées $\left(x,y\right)$
sur $\mathbb{C}^{2}$ dans lequel $\phi$ s'écrit sous la forme $\left(x,y\right)\mapsto\left(x,y+P\left(x\right)\right)$,
où $P\in\mathbb{C}\left[x\right]$. Il est connu depuis longtemps
\cite{Jun42} que tout automorphisme du plan peut s'écrire comme une
composition d'automorphismes affines et triangulaires. 

\begin{enavant} La conjecture naturelle, dite \emph{conjecture des
générateurs modérés}, serait de déclarer que tout automorphisme de
$\mathbb{C}^{3}$ s'obtient encore comme composition d'affines et
de triangulaires, un automorphisme $\phi$ de $\mathbb{C}^{3}$ étant
dit triangulaire s'il peut cette fois s'écrire sous la forme $\left(x,y,z\right)\mapsto\left(x,y+P\left(x\right),z+Q\left(x,y\right)\right)$,
où $P\in\mathbb{C}\left[x\right]$ et $Q\in\mathbb{C}\left[x,y\right]$,
dans un système de coordonnées $\left(x,y,z\right)$ bien choisi.
Un automorphisme admettant une décomposition en produit d'automorphismes
affines et triangulaires est dit \emph{modéré}. En 1972, Nagata \cite{Nag72}
a construit un candidat explicite destiné à infirmer la conjecture
des générateurs modérés, à savoir l'automorphisme $\sigma$ de $\mathbb{C}^{3}$
suivant :

\[
\left(x,y,z\right)\mapsto\sigma\left(x,y,z\right)=\left(x,y+x\left(xz+y^{2}\right),z-2y\left(xz+y^{2}\right)-x\left(xz+y^{2}\right)^{2}\right)\]
 Cet automorphisme est un exemple de ce que l'on appelle un automorphisme
de type exponentiel, c'est-à-dire un automorphisme dont le co-morphisme
est de la forme forme ${\rm exp}\left(\partial\right)$, où $\partial$
est une dérivation localement nilpotente. On voit en effet facilement
que l'automorphisme de Nagata correspond à la dérivation localement
nilpotente \[
\partial=\left(xz+y^{2}\right)\left(x\frac{\partial}{\partial y}-2y\frac{\partial}{\partial z}\right)\]
de $\mathbb{C}\left[x,y,z\right]$. Malgré de nombreux travaux consacrés
à cet automorphisme, il aura fallu attendre le début du vingtième
siècle pour que U. Umirbaev et I. Shestakov parviennent à démontrer
le fait suivant : 

\end{enavant}

\begin{thm}
\emph{(\cite{UmS04-2}, Corollaire 9)} L'automorphisme de Nagata n'est
pas modéré. 
\end{thm}
\begin{enavant} Il existe donc des automorphisme ``sauvages'' de
$\mathbb{C}^{3}$. Ce que les deux auteurs ont en réalité établi,
c'est que le fait qu'un automorphisme de $\mathbb{C}^{3}$ soit modéré
peut se décider de manière algorithmique. En particulier, leur critère
peut s'appliquer à une large classe d'automorphismes pour lesquels
l'algorithme en question garantit la non existence d'une décomposition
en produit d'automorphismes affines et triangulaires. L'automorphisme
de Nagata a cependant la propriété remarquable d'être ce que l'on
appelle stablement modéré. Cela signifie que si l'on étend cet automorphisme
en l'automorphisme $\tilde{\phi}=\phi\times{\rm id}$ de $\mathbb{C}^{4}$,
ce dernier se trouve être modéré. En effet, le co-morphisme de $\tilde{\phi}$
s'identifie à l'automorphisme exponentiel $\exp\left(\partial\right)$
de $\mathbb{C}\left[x,y,z,u\right]$, où l'on a considéré $\partial$
comme une dérivation localement nilpotente de $\mathbb{C}\left[x,y,z,u\right]$.
Or on a $\exp\left(\partial\right)=\tau^{-1}\exp\left(-\partial_{1}\right)\tau\exp\left(\partial_{1}\right)$,
où $\tau$ désigne l'automorphisme triangulaire de $\mathbb{C}\left[x,y,z,u\right]$
défini par $\tau\left(x,y,z,u\right)=\left(x,y,z,u+\left(xz+y^{2}\right)\right)$,
et où \[
\partial_{1}=u\left(x\frac{\partial}{\partial y}-2y\frac{\partial}{\partial z}\right)\]
 est une dérivation localement nilpotente triangulaire de $\mathbb{C}\left[x,y,z,u\right]$
induisant par conséquent un automorphisme exponentiel modéré. 

\end{enavant}

\subsection{Un nouvel exemple d'automorphisme sauvage stablement modéré}

\indent\newline\noindent L'automorphisme de Nagata apparaît comme
l'automorphisme exponentiel associé à la dérivation localement nilpotente
\[
\partial=\left(xz+y^{2}\right)\left(x\frac{\partial}{\partial y}-2y\frac{\partial}{\partial z}\right),\]
 obtenue à partir de la dérivation localement nilpotente triangulaire
${\displaystyle x\frac{\partial}{\partial y}-2y\frac{\partial}{\partial z}}$
annulant le polynôme $xz+y^{2}$. La surface d'équation $xz+y^{2}=0$
dans $\mathbb{C}^{3}$ constitue un cas particulier de surface du
type $S_{Q,n}$. Plus généralement, toute surface $S_{Q,n}$, définie
donc par une équation de la forme $x^{n}z-Q\left(x,y\right)=0$, est
annulée par la dérivation localement nilpotente triangulaire \[
\partial_{Q,n}=x^{n}\frac{\partial}{\partial y}+\frac{\partial Q\left(x,y\right)}{\partial y}\frac{\partial}{\partial z}\]
 de $\mathbb{C}\left[x,y,z\right]$. En se basant sur un critère établi
par U. Umirbaev et I. Shestakov \cite{UmS04-2}, on peut alors montrer
que pour tout entier $m\geq1$, l'automorphisme exponentiel associé
à une dérivation localement nilpotente du type $\left(x^{n}z-Q\left(x,y\right)\right)^{m}\partial_{Q,n}$
est sauvage. 

\begin{enavant} Les automorphismes sauvages construits à partir des
surfaces du type $S_{Q,n}$ par le procédé ci-dessus jouissent de
propriétés en tout point analogues à celles de l'automorphisme de
Nagata. En particulier, leur lieu fixe est supporté par la réunion
d'un certain nombre de droites dans $\mathbb{C}^{3}$ et de la surface
$S_{Q,n}$ correspondante. Nous allons construire ci-dessous à partir
de la surface $S\subset\mathbb{C}^{3}$ d'équation $P=xy-\left(f+1\right)\left(f+4\right)=0$,
où $f=xz-y^{2}\in\mathbb{C}\left[x,y,z\right]$ un automorphisme sauvage
stablement modéré de $\mathbb{C}^{3}$ qui ne sera conjugué à aucun
des automorphismes associé aux surfaces de type $S_{Q,n}$. 

\end{enavant}

\begin{enavant} Considérons pour cela la dérivation localement nilpotente
\[
\partial=x\left(2f+5\right)\frac{\partial}{\partial y}+\left(x+2y\left(2f+5\right)\right)\frac{\partial}{\partial z}=-\frac{\partial P}{\partial z}\frac{\partial}{\partial y}+\frac{\partial P}{\partial y}\frac{\partial}{\partial z}\]
 de $\mathbb{C}\left[x,y,z\right]$ définie dans le Théorème \ref{thm:Th_Principal}
ci-dessus. Bien que $\partial$ ne soit pas triangulaire dans le système
de coordonnées $\left(x,y,z\right)$ considéré, il résulte néanmoins
d'un critère dû à D. Daigle \cite{Dai96}, que $\partial$ est triangulable,
au sens qu'elle devient triangulaire dans un système de coordonnées
bien choisi. Par conséquent l'automorphisme exponentiel $\exp\left(\partial\right)$
est modéré. Par construction, la dérivation ci-dessus annule l'équation
$P=0$ de la surface $S\subset\mathbb{C}^{3}$. La dérivation $P\partial$
de $\mathbb{C}\left[x,y,z\right]$ est donc encore localement nilpotente,
mais n'est plus triangulable :

\end{enavant}

\begin{lem}
La dérivation localement nilpotente \[
P\partial=\left(xy-\left(f+1\right)\left(f+4\right)\right)\left(x\left(2f+5\right)\frac{\partial}{\partial y}+\left(x+2y\left(2f+5\right)\right)\frac{\partial}{\partial z}\right),\qquad\textrm{où }f=xz-y^{2}\]
 de $\mathbb{C}\left[x,y,z\right]$ n'est pas triangulable.
\end{lem}
\begin{proof}
Un calcul élémentaire (cf. Lemme 9.5.12 dans \cite{VdE00}) implique
que si $P\partial$ était triangulable, alors toute composante irréductible
du lieu fixe de l'automorphisme exponentiel associé serait un cylindre,
i.e. un schéma de la forme $Z\times\mathbb{C}$ pour un certain schéma
$Z$. Il résulte de \ref{txt:Lieu_fixe_de_partial} ci-dessus et de
la construction de la dérivation $P\partial$ que le lieu fixe dans
$\mathbb{C}^{3}$ de l'automorphisme associé à $\exp\left(P\partial\right)$
est la réunion des droites $L_{\pm}=\left\{ x=0;\, y=\pm\sqrt{5/2}\right\} $
et de la surface $S$ du Théorème \ref{thm:Th_Principal}. Or il n'existe
pas de courbe affine $C$ telle que $S$ soit isomorphe au cylindre
$C\times\mathbb{C}$. En effet, sinon $C$ devrait être une courbe
affine rationnelle puisque $S$ est une surface rationnelle. Le fait
que $H_{1}\left(S,\mathbb{Z}\right)=0$ impliquerait alors que $C$
est isomorphe à la droite affine complexe, donc que $S$ est isomorphe
au plan affine, ce qui est absurde. 
\end{proof}
\noindent Puisque la dérivation localement nilpotente $P\partial$
de $\mathbb{C}\left[x,y,z\right]$ n'est pas triangulable, l'automorphisme
exponentiel $\exp\left(P\partial\right)$ constitue un candidat d'automorphisme
sauvage. Il se trouve que c'est effectivement le cas, comme le montre
le résultat suivant. 

\begin{prop}
L'automorphisme de $\mathbb{C}^{3}$ suivant \[
\phi:\mathbb{C}^{3}\rightarrow\mathbb{C}^{3},\quad\left(\begin{array}{l}
x\\
y\\
z\end{array}\right)\mapsto\left(\begin{array}{l}
x\\
y+x\left(2f+5\right)P+x^{3}P^{2}\\
z+2P\left(2f+5\right)y+xP^{2}\left(2xy+2f+5\right)^{2}+2x^{3}P^{3}\left(2f+5\right)+x^{5}P^{4}\end{array}\right)\]
 est sauvage mais stablement modéré. 
\end{prop}
\begin{proof}
Nous laissons le soin au lecteur de vérifier que $\phi$ est simplement
l'automorphisme de $\mathbb{C}^{3}$ correspondant à l'automorphisme
$\exp\left(P\partial\right)$ de $\mathbb{C}\left[x,y,z\right]$.
Puisque $\partial\left(x\right)=0$, $\exp\left(P\partial\right)$
est un $\mathbb{C}\left[x\right]$-automorphisme de $\mathbb{C}\left[x,y,z\right]$.
D'après le Corollaire 10 dans \cite{UmS04-2}, $\exp\left(P\partial\right)$
est un automorphisme sauvage de $\mathbb{C}\left[x,y,z\right]$ si
et seulement si il l'est lorsqu'on le considère comme un $\mathbb{C}\left[x\right]$-automorphisme.
Cela revient donc à montrer que $\exp\left(P\partial\right)$ n'admet
pas de décomposition en produit de $\mathbb{C}\left[x\right]$-automorphismes
affines et triangulaires de $\mathbb{C}\left[x,y,z\right]$. Soient
$d_{2}=8$ et $d_{3}=16$ les degrés totaux en les variables $y$
et $z$ des composantes non triviales $\phi_{2}$ et $\phi_{3}$ de
$\phi$ et soient $F_{2}=x^{3}y^{8}$ et $F_{3}=x^{5}y^{16}$ les
composantes homogènes de degrés $d_{2}$ et $d_{3}$ des polynômes
$\phi_{2}$ et $\phi_{3}$. Si $\exp\left(P\partial\right)$ était
un $\mathbb{C}\left[x\right]$-automorphisme modéré de $\mathbb{C}\left[x,y,z\right]$
alors, en vertue du Corollaire 5.1.6 dans \cite{VdE00}, il existerait
$c\in\mathbb{C}\left[x\right]$ tel que $x^{5}y^{16}=F_{3}=cF_{1}^{2}=cx^{6}y^{16}$,
ce qui est clairement impossible. Ainsi $\phi$ est un automorphisme
sauvage de $\mathbb{C}^{3}$. Pour conclure la preuve, il suffit d'établir
que l'automorphisme $\tilde{\phi}=\phi\times{\rm id}$ de $\mathbb{C}^{4}$
est modéré. L'argument que nous allons utiliser repose sur une méthode
générale due à M Smith \cite{Smii89}. Soit $\tilde{\partial}$ la
dérivation localement nilpotente triangulable de $\mathbb{C}\left[x,y,z,u\right]$
obtenue à partir de la dérivation $\partial$ du Théorème \ref{thm:Th_Principal}
en posant $\tilde{\partial}\left(u\right)=0$. Puisque $u\in{\rm Ker}\left(\tilde{\partial}\right)$,
la dérivation $u\tilde{\partial}$ est encore localement nilpotente.
On observe alors que $\left(\exp\left(P\partial\right)\otimes{\rm id}\right)=\exp\left(P\tilde{\partial}\right)=\tau^{-1}\exp\left(-u\tilde{\partial}\right)\tau\exp\left(u\tilde{\partial}\right)$
où $\tau$ désigne l'automorphisme triangulaire de $\mathbb{C}^{4}$
défini par $\left(x,y,z,u\right)\mapsto\left(x,y,z,u+P\left(x,y,z\right)\right)$.
Il suffit par conséquent de vérifier que les automorphismes $\exp\left(\pm u\tilde{\partial}\right)$
de $\mathbb{C}\left[x,y,z,u\right]$ sont modérés. Or puisque $\partial$
est triangulable, la dérivation localement nilpotente $u\tilde{\partial}$
de $\mathbb{C}\left[x,y,z,u\right]$ l'est aussi. Ainsi $\exp\left(u\tilde{\partial}\right)$
et son inverse $\exp\left(-u\tilde{\partial}\right)$ sont bien modérés.
\end{proof}
\begin{enavant} Le lecteur se convaincra sans peine que des arguments
analogues à ceux utilisés ci-dessus montreraient que pour tout couple
de polynômes $Q_{1},Q_{2}\in\mathbb{C}\left[x,y\right]$ et tout couple
d'entiers $\left(n_{1},n_{2}\right)\in\mathbb{Z}_{\geq1}^{2}$, la
surface $S_{\left(Q_{1},Q_{2}\right),\left(n_{1},n_{2}\right)}$ de
$\mathbb{C}^{3}$ définie par l'équation $P=x^{n_{2}}y-Q_{2}\left(x,f\right)=0$,
où $f=x^{n_{1}}z-Q_{1}\left(x,y\right)$, est invariante par une dérivation
localement nilpotente $\partial$ de $\mathbb{C}\left[x,y,z\right]$
pour laquelle l'automorphisme exponentiel $\exp\left(P\partial\right)$
associé sera sauvage mais stablement modéré. Nous renvoyons le lecteur
intéressé à un futur article \cite{DubPrep} dans lequel les surfaces
$S_{\left(Q_{1},Q_{2}\right),\left(n_{1},n_{2}\right)}$ de $\mathbb{C}^{3}$
introduites ci-dessus seront étudiées en détails. 

\end{enavant}

\bibliographystyle{amsplain}

\begin{thebibliography}{10}

\bibitem{BML01}
T.~Bandman and L.~Makar-Limanov, \emph{Affine surfaces with
  ${AK}\left({S}\right)=\mathbb{C}$}, Michigan J. Math. \textbf{49} (2001),
  567--582.

\bibitem{Ber83}
J.~Bertin, \emph{Pinceaux de droites et automorphismes des surfaces affines},
  J. reine angew. Math. \textbf{341} (1983), 32--53.

\bibitem{Bon02}
P.~Bonnet, \emph{Surjectivity of quotient maps for algebraic $(\mathbb{C}
  ,+)$-actions and polynomial maps with contractible fibers}, Transformation
  {G}roups \textbf{7} (2002), no.~1, 3--14.

\bibitem{Dai96}
D.~Daigle, \emph{A necessary and sufficient condition for triangulability of
  derivations of $k\left[x,y,z\right]$}, J. of {P}ure and {A}pplied {A}lgebra
  \textbf{113} (1996), 297--305.

\bibitem{DaRus02}
D.~Daigle and P.~Russel, \emph{On log $\mathbb{Q}$-homology planes and weighted
  projective planes}, To appear in {C}anadian {J}. {M}ath.

\bibitem{Dan89}
W.~Danielewski, \emph{On a cancellation problem and automorphism groups of
  affine algebraic varieties}, Preprint {W}arsaw, 1989.

\bibitem{DubPrep}
A.~Dubouloz, \emph{Makar-limanov surfaces in $\mathbb{C}^3$}, In preparation.

\bibitem{Dub02}
\bysame, \emph{Completions of normal affine surfaces with a trivial
  {M}akar-{L}imanov invariant}, Michigan J. Math. \textbf{52} (2004), no.~2,
  289--308.

\bibitem{DubG03}
\bysame, \emph{{D}anielewski-{F}ieseler {S}urfaces}, {T}ransformation {G}roups
  \textbf{10} (2005), no.~2, 139--162.

\bibitem{DubP04}
A.~Dubouloz and P.M. Poloni, \emph{On a class of {D}anielewski surfaces in
  affine $3$-space}, In preparation.

\bibitem{Fie94}
K.H. Fieseler, \emph{On complex affine surfaces with $\mathbb{C}_+$-actions},
  Comment. Math. Helvetici \textbf{69} (1994), 5--27.

\bibitem{FrMo02}
G.~Freudenburg and L.~Moser-Jauslin, \emph{Embeddings of {D}anielewski
  surfaces}, Math. Z. \textbf{245} (2003), no.~4, 823--834.

\bibitem{Giz71}
M.H. Gizatullin, \emph{Quasihomogeneous affine surfaces}, Math. USSR Izvestiya
  \textbf{5} (1971), 1057--1081.

\bibitem{GrHar78}
P~Griffiths and J.~Harris, \emph{Principles of {A}lgebraic {G}eometry}, John
  {W}iley, New York, 1978.

\bibitem{SGA1}
A.~Grothendieck~et al., \emph{{SGA1}. {R}ev{ê}tements {é}tales et {G}roupe
  {F}ondamental}, Lecture {N}otes in {M}ath., vol. 224, Springer-Verlag, 1971.

\bibitem{GuMiy05}
R.~V. Gurjar and M.~Miyanishi, \emph{Automorphisms of affine surfaces with
  $\mathbb{A}^1$-fibrations}, Michigan Math. J. \textbf{53} (2005), no.~1,
  33--55.

\bibitem{Jun42}
H.W.E. Jung, \emph{{Ü}ber ganze birationale {T}ransformationen der {E}bene}, J.
  {R}eine {A}ngew. {M}ath. \textbf{184} (1942), 161--174.

\bibitem{KaZa99}
S.~Kaliman and M.~Zaidenberg, \emph{Affine modifications and affine
  hypersurfaces with a very transitive automorphism group}, Transformation
  Groups \textbf{4} (1999), 53--95.

\bibitem{ML01}
L.~Makar-Limanov, \emph{On the group of automorphisms of a surface
  $x^ny=p\left(z\right)$}, Israel J. Math. \textbf{121} (2001), 113--123.

\bibitem{Miy85}
M.~Miyanishi, \emph{Normal affine subalgebras of a polynomial ring}, Algebraic
  and Topological Theories. To The Memory of Dr. Takehito Miyata, Kinokuniya,
  Tokyo, 1985, pp.~37--51.

\bibitem{MoP05}
L.~Moser-Jauslin and P-M. Poloni, \emph{Emebddings of a family of {D}anielewski
  surfaces and certain $\mathbb{C}_+$-actions on $\mathbb{C}^3$},  (2005),
  Preprint {U}niversit\'e de {B}ourgogne.

\bibitem{Nag72}
M.~Nagata, \emph{On the automorphism group of $k[x,y]$}, {K}yoto {U}niv.
  {L}ectures in {M}ath., vol.~5, {K}yoto {U}niversity , {K}inokuniya-{T}okyo,
  1972.

\bibitem{UmS04-2}
I.~P. Shestakov and U.~U. Umirbaev, \emph{The tame and the wild automorphisms
  of polynomial rings in three variable}, J. Amer. Math. Soc. \textbf{17}
  (2004), no.~197--227.

\bibitem{Smii89}
M.K. Smith, \emph{Stably tame automorphisms}, J. of {P}ure and {A}pplied
  {A}lgebra \textbf{58} (1989), 209--212.

\bibitem{VdE00}
A.~van~den Essen, \emph{Polynomial automorphisms and the {J}acobian
  conjecture}, Progress in {M}athematics, vol. 190, Birkh\"auser {V}erlag,
  Basel, 2000.

\end{thebibliography}

\end{document}